\documentclass{paper}

\usepackage{hyperref}
\usepackage{graphicx}
\usepackage{amsmath}
\usepackage{xcolor}
\usepackage{amssymb}
\usepackage{amsfonts}
\usepackage{marginnote}
\usepackage{fancybox}
\usepackage{color}
\usepackage{bm}
\usepackage{cancel}
\usepackage{wrapfig}
\usepackage{float}

\newtheorem{lemma}{Lemma}[section]
\newtheorem{theorem}{Theorem}[section]
\newtheorem{remark}{Remark}[section]
\newtheorem{corollary}{Corollary}[section]

\evensidemargin=.in

\def\gam{\gamma}

\newcommand{\mcO}{\mathcal{O}}

\newcommand{\mcD}{\mathcal{D}}
\newcommand{\mcG}{\mathcal{G}}

\newcommand{\mcL}{\mathcal{L}}

\newcommand{\mbR}{\mathbb{R}}
\newcommand{\mbRd}{{\mathbb{R}^d}}

\newcommand{\omg}{{\omega}}
\newcommand{\omgc}{{\eta_c}}
\newcommand{\omgi}{{\eta_D}}
\newcommand{\gamc}{{\Gamma_c}}
\newcommand{\gami}{{\Gamma_D}}

\newcommand{\omgn}{{\omega_n}}
\newcommand{\omgnp}{{\Omega_n}}
\newcommand{\omgl}{{\Omega_l}}
\newcommand{\omgb}{{\Omega_o}}

\newcommand{\omgw}{{\eta_n}}
\def\omgnw{{\omg_n\hspace{-.4mm}\cup\omgw}}
\newcommand{\womg}{{\eta}}

\def\omgp{{\Omega}}

\newcommand{\beq}{\begin{equation}}
\newcommand{\eeq}{\end{equation}}

\def \bphi{{\boldsymbol\Phi}}

\def \alphab{{\boldsymbol\alpha}}
\def \nub{{\boldsymbol \nu}}

\def \xb{\bm{x}}
\def \yb{\bm{y}}

\def \ul {u_l}
\def \un {u_n}

\def \vl {v_l}
\def \vn {v_n}
\def \zn {z_n}
\def \zl {z_l}
\def \ulh {{u_l^h}}
\def \unh {{u_n^h}}
\def \uloh {{u_l^{h0}}}
\def \unoh {{u_n^{h0}}}

\def \vlh {{v_l^h}}
\def \vnh {{v_n^h}}
\def \znh {{z_n^h}}
\def \zlh {{z_l^h}}
\def \thn {\theta_n}
\def \thl {\theta_l}

\def \thnh {\theta_n^h}
\def \thlh {\theta_l^h}
\def \thna {\theta_n^*}
\def \thla {\theta_l^*}
\def \wun {\widehat u_n}
\def \wul {\widehat u_l}
\def \unl {u^{*}}
\def \sgn {\sigma_n}
\def \sgl {\sigma_l}
\def \sgnh {{\sigma_n^h}}
\def \sglh {{\sigma_l^h}}

\def\mydate{\number\day\ {\ifcase\month \or January\or February\or
           March\or April\or May\or June\or July\or August\or
            September\or October\or November\or December\fi}
\number\year}

\title{\textsf{\textbf{Formulation, analysis and computation of an optimization-based local-to-nonlocal coupling method}}}

\author{
Marta D'Elia\thanks{Computational Science and Analysis, Sandia National Laboratories, Livermore, CA, 94550 ({\tt mdelia@sandia.gov}).} 
\and
Pavel Bochev\thanks{Center for Computing Research, Sandia National Laboratories, Albuquerque, NM, 87185 ({\tt pbboche@sandia.gov}).}}

\begin{document}
\maketitle

\begin{abstract}
We present an optimization-based coupling method for local and nonlocal continuum models. Our approach couches the coupling of the models into a control problem where the states are the solutions of the nonlocal and local equations, the objective is to minimize their mismatch on the overlap of the local and nonlocal problem domains, and the virtual controls are the nonlocal volume constraint and the local boundary condition. 
We present the method in the context of Local-to-Nonlocal diffusion coupling. Numerical examples illustrate the theoretical properties of the approach.
\end{abstract}
\vspace{-2ex}

\pagestyle{myheadings}
\thispagestyle{plain}
\markboth{OPTIMIZATION-BASED LOCAL-TO-NONLOCAL COUPLING}{OPTIMIZATION-BASED LOCAL-TO-NONLOCAL COUPLING}

\section{Introduction} \label{intro}

Nonlocal continuum theories such as peridynamics \cite{Silling_10_AAM},  physics-based nonlocal elasticity \cite{DiPaola_09_JE}, or nonlocal descriptions resulting from homogenization of nonlinear damage models \cite{Han_12_IJNME} can incorporate strong nonlocal effects due  to long-range forces at the mesoscale or microscale. 
As a result, for problems where these effects cannot be neglected, such descriptions are more accurate  than local Partial Differential Equations (PDEs) models. However, their computational cost is also significantly higher than that of PDEs.
Local-to-Nonlocal (LtN) coupling methods aim to combine the computational efficiency of PDEs with the accuracy of nonlocal models.
The need for LtN couplings is especially acute when the size of the computational domain is such that the nonlocal solution becomes prohibitively expensive to compute, yet the nonlocal model is required to accurately resolve small scale features such as crack tips or dislocations that can affect the global material behavior \cite{DElia_2020review}.

LtN couplings involve two fundamentally different mathematical descriptions of the same physical phenomena. The principal challenge is the stable and accurate merging of these descriptions into a physically consistent coupled formulation.
In this paper we address this challenge by couching the LtN coupling into an optimization problem. The objective is to minimize the  mismatch of the local and nonlocal solutions on the overlap of their respective subdomains, the constraints are the associated governing equations, and the controls are the virtual nonlocal volume constraint and the local boundary condition. We formulate and analyze this optimization-based LtN approach in the context of local and nonlocal diffusion models  \cite{Du_12_SIREV}.

Our coupling strategy differs fundamentally from other LtN approaches such as the extension of the Arlequin \cite{BenDhia_05_IJNME} method to LtN couplings \cite{Han_12_IJNME}, force-based couplings \cite{Seleson_13_CMS}, or the morphing approach \cite{Azdoud_13_IJSS,Lubineau_12_JMPS}. The first two schemes blend the energies or the forces of the two models over a dedicated ``gluing'' area, while the third one implements the coupling through a gradual change in the material properties characterizing the two models over a ``morphing'' region.   In either case, resulting LtN methods treat the coupling condition as a constraint, similar to classical domain decomposition methods.
In contrast, we treat this condition as an optimization objective, and keep the two models separate. This strategy brings about valuable theoretical and computational advantages.  For instance, the coupled problem passes a patch test by construction and its well-posedness typically follows from the well-posedness of the constraint equations. 
 
Our approach has its roots in non-standard  \emph{optimization-based domain decomposition} methods for PDEs  \cite{Discacciati_13_SIOPT,Du_01_SINUM,Du_00_SINUM,Gervasio_01_NM,Gunzburger_00_AMC,Gunzburger_00_SINUM,Gunzburger_99_CMA,Kuberry_13_CMAME}. It  has also been applied to the coupling of discrete atomistic and continuum models in \cite{Bochev_14_SINUM,Bochev_14_LSSC} and multiscale problems \cite{Abdulle_15}.
This paper continues the efforts in \cite{Bochev_14_INPROC}, which presented an initial optimization-based LtN  formulation, in 
\cite{Bochev_16b_CAMWA}, which focussed on specializing the formulation to mixed boundary conditions and mixed volume constraints and its practical demonstration using Sandia's agile software components toolkit, and in \cite{DElia_15Handbook}, which extended the formulation to peridynamics.
The main contributions of this paper include (i) rigorous analysis of the LtN coupling error, (ii) formal proof of the well-posedness of the discretized LtN formulation, and (iii) rigorous convergence analysis.

We have organized the paper as follows. Section \ref{notation} introduces notation, basic notions of nonlocal vector calculus and the relevant mathematical models. We present the optimization-based LtN method and prove its well-posedness in Section \ref{optimization} and study its error in Section \ref{error-analysis}. Section \ref{finite-dim} focusses on the discrete LtN formulation, its well-posedness and numerical analysis. A collection of numerical examples in Section \ref{numerical-tests} illustrates the theoretical properties of the method using a simple one-dimensional setting.

\section{Preliminaries} \label{notation}
Let $\omega$ be a bounded open domain in $\mbRd$, $d=2,3$, with Lipschitz-continuous boundary $\partial\omega$. We use the standard notation $H^s(\omega)$ for a Sobolev space of order $s$ with norm and inner product $\|\cdot\|_{s,\omega}$ and $(\cdot,\cdot)_{s,\omega}$, respectively. As usual, $H^0(\omega):=L^2(\omega)$ is the space of all square integrable functions on $\omega$.  The subset of all functions in $H^1(\omega)$ that vanish on $\partial\omega$ is $H^1_0(\omg) := \{u\in H^1(\omega): \; u|_{\partial\omega}=0 \}$.

\begin{figure}
\centering
\includegraphics[width=0.3\textwidth]{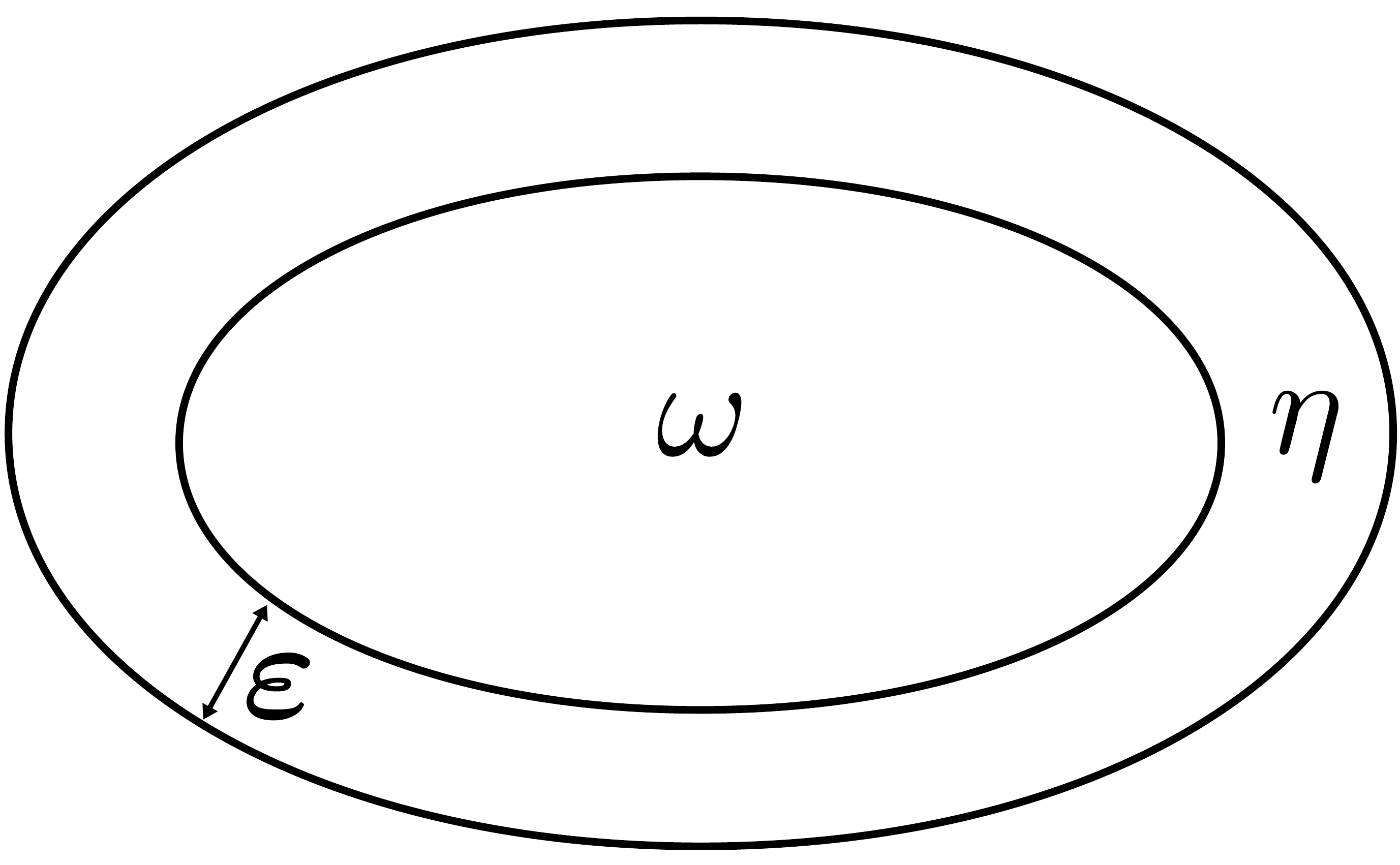}
\vspace{-1ex}
\caption{Two-dimensional domain $\omg$ and interaction domain $\womg$ with interaction radius $\varepsilon$.}
\vspace{-2ex}
\label{full-domain}
\end{figure}

The nonlocal model in this paper requires nonlocal vector calculus operators \cite[\S3.2]{Du_12_SIREV} acting on functions
$u(\xb)\colon \mbRd \to \mbR$ and $\nub(\xb,\yb)\colon \mbRd\times\mbRd\to \mbRd$.
Let $\gamma (\xb, \yb )\colon\mbRd\times\mbRd\to\mbR$ and $\alphab(\xb,\yb) \colon \mbRd\times\mbRd\to \mbRd$ be a non-negative symmetric kernel and an antisymmetric function, respectively, i.e.,  $\gamma(\xb, \yb )=\gamma(\yb, \xb )\ge0$ and $\alphab(\yb,\xb)=-\alphab(\xb,\yb)$. 
The nonlocal \emph{diffusion}\footnote{More general representations of $\mcL$, associated with non-symmetric and not necessarily positive kernel functions exist. Such nonlocal operators may define models  for non-symmetric diffusion phenomena such as non-symmetric jump processes \cite{DElia_17CMAM}.} of $u$ is an operator $\mcL(u)\colon \mbRd \to \mbR$ defined by  
\begin{displaymath}
\mcL u(\xb) := 2\int_\mbRd \big(u(\yb)-u(\xb)\big) \, \gamma (\xb, \yb )\,d\yb \qquad  \xb \in \mbRd,
\end{displaymath}
and its nonlocal \emph{gradient} is a mapping $\mcG(u)\colon  \mbRd\times\mbRd\to\mbRd$ given by
\begin{subequations}\label{eq:ndivgrad}
\begin{equation}\label{ngra}
\mcG u(\xb,\yb) := \big(u(\yb)-u(\xb)\big)  \alphab(\xb,\yb) \qquad \xb,\yb\in\mbRd.
\end{equation}
Finally, the nonlocal \emph{divergence} of $\nub(\xb,\yb)$ is a mapping $\mcD(\nub)\colon \mbRd \to \mbR$ defined by\footnote{The paper \cite{Du_13_MMMAS} shows that the adjoint $\mcD^*=-\mcG$.}
\begin{equation}\label{ndiv}
\mcD\nub(\xb) := \int_{\mbRd} \big(\nub(\xb,\yb)+\nub(\yb,\xb)\big)\cdot\alphab(\xb,\yb)\,d\yb\qquad \xb\in\mbRd.
\end{equation}
\end{subequations}
Furthermore, given a second-order symmetric tensor $\bphi(\xb,\yb)=\bphi(\yb,\xb)$, equations \eqref{eq:ndivgrad} imply that
\begin{displaymath}
\mcD\big(\bphi \mcG u)(\xb)  = 2\int_{\mbRd}\big(u(\yb)-u(\xb)\big) \big( \alphab(\xb,\yb)\cdot \bphi(\xb,\yb)\alphab(\xb,\yb)\big) \,d\yb.
\end{displaymath}
Thus, with the identification $\gamma(\xb,\yb):=\alphab(\xb,\yb)\cdot \bphi(\xb,\yb) \alphab(\xb,\yb)$
the operator $\mcL$ is a composition of the nonlocal divergence and gradient operators: $\mcL u = \mcD\big(\bphi \mcG u)$.
We define the \emph{interaction domain} of an open bounded region $\omega\in\mbRd$ as
\begin{displaymath}
\womg = \{\yb\in\mbRd\setminus\omega: \; \gamma(\xb,\yb)\neq 0\},
\end{displaymath}
for $\xb\in\omega$ and set $\Omega =\omega\cup\womg$. In this paper we consider kernels $\gamma$ such that for $\xb\in\omg$
\begin{equation}\label{eq:gamma-conds}
\left\{\begin{aligned}
\gamma(\xb,\yb)\geq  0 \quad &\forall\, \yb\in B_\varepsilon(\xb) \\
\gamma(\xb,\yb)    = 0 \quad &\forall\, \yb\in \mbRd \setminus B_\varepsilon(\xb),
\end{aligned}\right.
\end{equation}
where $B_\varepsilon(\xb) = \{\yb\in \mbRd: \; |\xb-\yb|\leq\varepsilon\}$. Kernels that satisfy \eqref{eq:gamma-conds} are referred to as localized kernels with {\it interaction radius} $\varepsilon$. It is easy to see that for such kernels the interaction domain is a layer of thickness $\varepsilon$ that surrounds $\omega$, i.e.
\begin{displaymath}
\womg = \{ \yb\in \mbRd\setminus\omega: \; \inf_{\xb\in\omega}|\yb-\xb|\leq\varepsilon\};
\end{displaymath}
see Fig. \ref{full-domain} for a two-dimensional example.
For a symmetric positive definite tensor $\bphi$ we respectively define the nonlocal energy semi-norm, nonlocal energy space, and nonlocal volume-constrained energy space by
\begin{subequations}
\begin{align}
&|||v|||_{\omgp}^2 := \frac12\int_{\omgp}\int_{{\omgp}} \mcG 
                      v\cdot(\bphi\mcG v)\,d\yb \, d\xb \label{energynorm}\\
&V({\omgp})        := \left\{ v  \in L^2(\omgp) \,\,:\,\, 
                      |||v|||_{\omgp} < \infty \right\} \label{vspace}\\
&V_{\tau}({\omgp}) := \left\{v\in V({\omgp}) \,\,:\,\, v=0\;{\rm on}\;\tau\right\}, 
                      \;\; \hbox{ for $\tau \subseteq \womg$.} \label{vcspace}
\end{align}
\end{subequations}
We also define the volume-trace space $\widetilde V(\womg):=\{v|_\womg: \,v\in V({\omgp})\}$, and an associated norm
\begin{equation}\label{eq:trace-norm} 
\|\sigma\|_{\widetilde V(\womg)}:=
\inf\limits_{v\in V({\omgp}), v|_\womg=\sigma} \, |||v |||_{\omgp}.
\end{equation}
We refer to  \cite{Du_12_SIREV,Du_13_MMMAS} for further information about the nonlocal vector calculus. 

In order to avoid technicalities not germane to the coupling scheme, in this paper we consider integrable kernels. Examples of applications modeled by the latter can be found in \cite{Aksoylu:2010,Aksoylu:2011,Andreu:2010}. Specifically, we assume that there exists positive constants $\gam_0$ and $\gam_2$ such that
\begin{equation}\label{eq:integrable-kernel}
\gam_0\leq \int_{\omgp\cap B_\varepsilon(\xb)} \gam(\xb,\yb)\,d\yb 
\qquad {\rm and} \qquad
\int_\omgp \gam(\xb,\yb)^2d\yb \leq \gam_2^2,
\end{equation}
for all $\xb\in\omgp$.
Note that this also implies that there exists a positive constant $\gamma_1$ such that for all $\xb\in\omgp$
\begin{equation}\label{eq:gamma1}
\int_\omgp \gam(\xb,\yb)\,d\yb \leq \gam_1.
\end{equation}
In \cite[\S4.2]{Du_12_SIREV} this class of kernels (referred to as Case 2) is rigorously analyzed; we report below an important result, useful throughout the paper.
\begin{lemma}\label{lem:L2-equivalence}
Let the function $\gam$ satisfy \eqref{eq:gamma-conds} and \eqref{eq:integrable-kernel}, then, there exist positive constants $C_{pn}$ and $C^*$ such that 
\begin{equation}\label{eq:L2equivalence}
\dfrac{1}{C_{pn}} \|u\|_{0,\omgp} \leq |||u|||_\omgp \leq C^* \|u\|_{0,\omgp}
\quad \forall\,u\in V_\tau(\omgp).
\end{equation}
Furthermore, the energy space $V_\tau(\omgp)$ is equivalent to $L_\tau^2(\omgp)=\{v\in L^2(\omgp): v|_\tau=0\}$.
\end{lemma}

\smallskip The latter is a combination of results in Lemmas 4.6 and 4.7 and Corollary 4.8 in \cite[\S4.3.2]{Du_12_SIREV}. Note that the lower bound in \eqref{eq:L2equivalence} represents a nonlocal Poincar\'e inequality.
Even though not included in the analysis, singular kernels appear in applications such as peridynamics; numerical results, included in the paper, suggest that the coupling scheme can handle such kernels without difficulties. However, their analysis is beyond the scope of this paper.

\subsection{Local-to-Nonlocal coupling setting} \label{nvcp}
Consider a bounded open region $\omega\subset\mbRd$ with interaction domain $\womg$. Given $f_n\in L^2(\omg)$ and $\sigma_n\in \widetilde V(\womg)$ we assume that  the volume-constrained\footnote{The volume constraint in \eqref{nonlocal-problem} is the nonlocal analogue of a Dirichlet boundary condition, in which the closed region $\womg$ plays the role of a``boundary''.} nonlocal  diffusion equation
\begin{equation} \label{nonlocal-problem}
\left\{\begin{array}{rlll}
-\mcL u_n & = & f_n & \xb\in \omg \\
u_n & = & \sigma_n & \xb\in\womg,
\end{array}\right.
\end{equation}
provides an accurate description of the relevant physical processes in $\omgp=\omg\cup\womg $. Let $\Gamma=\partial\omgp$, we assume  that the local diffusion model given by the Poisson equation
\begin{equation} \label{local-problem}
\left\{\begin{array}{rlll}
-\Delta u_l & = & f_l &  \xb\in \omgp \\
u_l  & = & \sigma_l &  \xb\in \Gamma,
\end{array}\right.
\end{equation} 
with suitable boundary data $\sigma_l\in H^\frac12(\Gamma)$ and forcing term $f_l\in L^2(\omgp)$  is a good approximation of  \eqref{nonlocal-problem} whenever the latter has sufficiently ``nice'' solutions. In this work we define $f_l$ to be an extension of $f_n$ by $0$  in $\womg$, specifically,
\begin{equation}\label{fl-def}
f_l=\left\{
\begin{aligned}
f_n & \quad \xb\in\omg\\
0   & \quad \xb\in\womg.
\end{aligned}\right.
\end{equation}
For a symmetric positive definite $\bphi$ standard arguments of variational theory show that the weak form\footnote{Multiplication of \eqref{nonlocal-problem} by a test function $\zn\in V_\womg(\omgp)$, integration over $\omg$ and application of the first nonlocal Green's identity \cite{Du_13_MMMAS} yield the weak form \eqref{nonlocal-weak} of the nonlocal problem.}
\begin{equation}\label{nonlocal-weak}
\int_{\omgp}\int_{\omgp} \mcG\un \cdot (\bphi\mcG\zn)\, d\yb \, d\xb = \int_\omg f_n \zn \, d\xb 
\qquad \forall\,\zn\in V_\womg({\omgp}) \,
\end{equation}
of \eqref{nonlocal-problem} is well-posed \cite{Du_12_SIREV}, i.e., \eqref{nonlocal-weak} has a unique solution such that 
\begin{equation}\label{cont-data-nonlocal}
|||\un|||_{\omgp} \leq K_n (\|f_n\|_{0,\omg}+\|\sgn\|_{\widetilde V(\womg)})
\end{equation}
for some positive constant $K_n$. In this work, for simplicity and without loss of generality, we set $\bphi={\bf I}$.

Although \eqref{nonlocal-weak} and the nonlocal calculus \cite{Du_12_SIREV} enable formulation and analysis of finite elements for \eqref{nonlocal-problem}, which parallel those for  the Poisson equation \eqref{local-problem}, resulting methods may be computationally intractable for large domains. The root cause for this is that long-range interactions increase the density of the resulting algebraic system making it more expensive to assemble and solve.

\section{Optimization-based LtN formulation}\label{optimization}
\begin{figure}
\centering
{\includegraphics[width=0.55\textwidth]{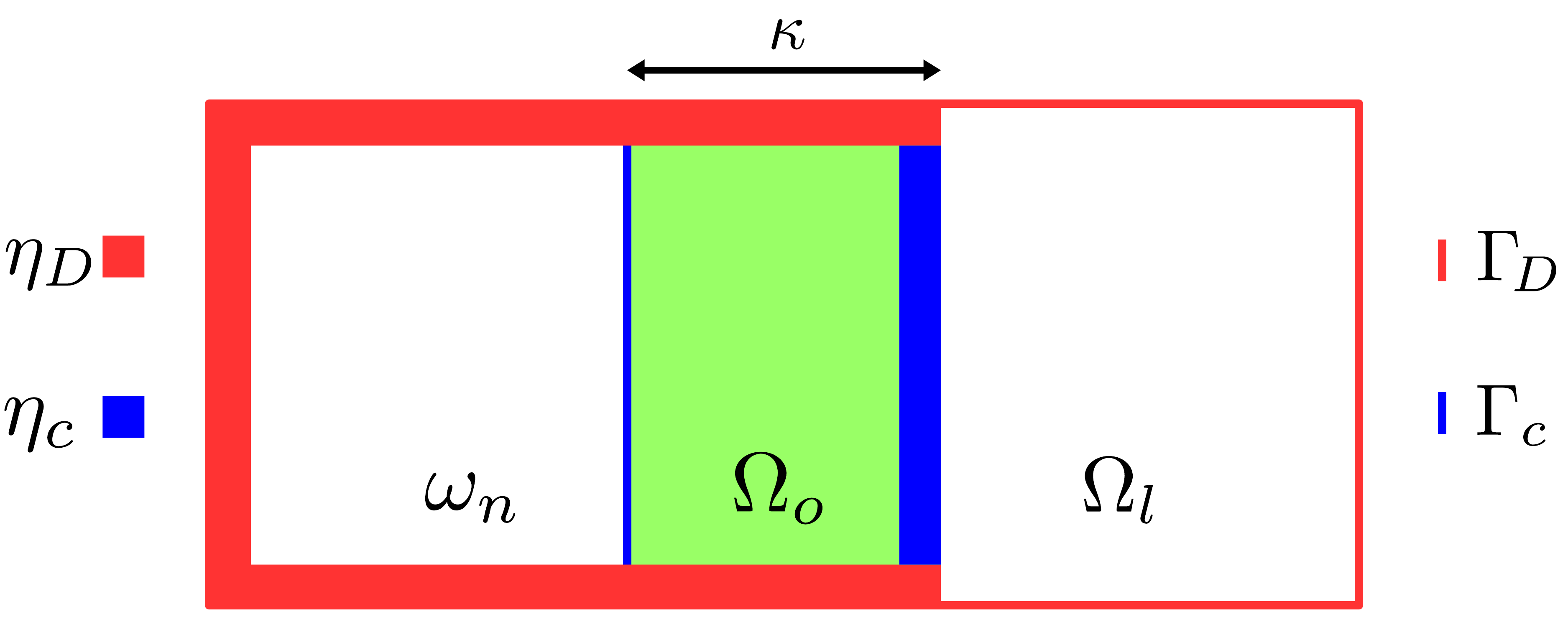}} 
\vspace{-2ex}
\caption{An example LtN domain configuration in two-dimensions.}
\vspace{-2ex}
\label{domains}
\end{figure}

For clarity we consider \eqref{nonlocal-problem} and \eqref{local-problem} with homogeneous Dirichlet conditions on $\womg$ and $\Gamma$, respectively.
To describe the approach it suffices to examine a coupling scenario where these problems operate on two overlapping subdomains of $\omgp$. 
Thus we consider partitioning of $\omgp$ into  a nonlocal subdomain $\omgn$ with interaction volume $\omgw$ and a local subdomain $\omgl$, with boundary $\Gamma_l$, such that $\omgnp:=\omgnw\subset\omgp$, $\Omega=\omgnp\cup\omgl$ and $\omgb = \omgnp\cap\omgl\neq\emptyset$; see Fig. \ref{domains}.
Let $\omgi=\womg\cap\omgw$, $\omgc=\omgw\setminus\omgi$, $\gami=\Gamma\cap\Gamma_l$, and $\gamc=\Gamma_l\setminus\gami$; see Fig. \ref{domains} and Appendix \ref{app:notation} for a summary of notation and definitions. Restrictions of \eqref{nonlocal-problem} and \eqref{local-problem} to $\omgn$ and $\omgl$ are given by
\begin{equation}\label{eq:subprob}
\left\{\begin{array}{rlll}
-\mcL\un & = & f_n & \xb\in\omgn\\
\un & = & \thn & \xb\in\omgc\\
\un & = & 0 & \xb\in\omgi
\end{array}\right.
\quad\mbox{and}\quad
\left\{\begin{array}{rlll}
-\Delta \ul & = & f_l & \xb\in\omgl\\ 
\ul & = & \thl & \xb\in\gamc\\ 
\ul & = & 0 & \xb\in\gami, 
\end{array}\right.
\end{equation}
respectively, where $\thn\in \Theta_n=\{\vn|_\omgc: \vn\in V_\omgi(\omgnp)\}$ and
$\thl\in\Theta_l=\{\vl|_\gamc:\vl\in H^1_\gami(\omgl)\}$ 
are an undetermined Dirichlet volume constraint and an undetermined Dirichlet boundary condition, respectively. 
The following constrained optimization problem
\begin{equation}\label{minimization}
\begin{array}{c}
\displaystyle\min\limits_{\un,\ul,\thn,\thl}  \displaystyle J(\un,\ul)
\;\;
\mbox{subject to \eqref{eq:subprob},} 
\quad \hbox{where} \;\;
J(\un,\ul) =\frac12 \| u_n-u_l \|^2_{0,\omgb}
\end{array}
\end{equation}
defines the optimization-based LtN coupling. In this formulation the subdomain problems  \eqref{eq:subprob} are the optimization constraints, $\un$ and $\ul$ are the states and $\thn$ and $\thl$ are the controls. We equip the control space $\Theta_n\times\Theta_l$ with the norm
\begin{equation}\label{eq:control-norm}
\|(\sigma_n,\sigma_l)\|^2_{\Theta_n\times\Theta_l} = 
\|\sgn\|^2_{0,\eta_c} + \|\sgl\|^2_{\frac12,\gamc}.
\end{equation}
In contrast to blending, \eqref{minimization} is an example of a divide-and-conquer strategy as the local and nonlocal problems operate independently in $\omgnp$ and $\omgl$.

Given an optimal solution $(\un^*,\ul^*, \thna,\thla)\in  V_\omgi(\omgnp)\times H^1_\gami(\omgl)\times\Theta_n\times\Theta_l$ of  \eqref{minimization} we define the LtN solution $\unl \in L^2(\omgp)$ by splicing together the optimal states:
\begin{equation}\label{unl-definition}
\unl = \left\{
\begin{array}{ll}
\un^* & \xb\in\omgnp\\ 
\ul^* & \xb\in\omgl\setminus\omgb.
\end{array}\right.
\end{equation}
We note that there are several ways to define the LtN solution; since our ultimate goal is to approximate a globally nonlocal solution, we set $u^*$ equal to the optimal nonlocal solution over the whole nonlocal domain.

The next section verifies that \eqref{minimization} is well-posed.

\subsection{Well-posedeness}\label{well-posedness}
We show that for any pair of controls  subproblems \eqref{eq:subprob} have unique solutions $\un(\thn)$ and $\ul(\thl)$, respectively. Elimination of the states from \eqref{minimization} yields the equivalent reduced space form of this problem in terms of $\thn$ and $\thl$ only:
\begin{equation}\label{reduced-minimization}
\displaystyle\min\limits_{\thn,\thl} J(\thn,\thl) \quad \hbox{where} \quad
J(\thn,\thl) = \frac12 \| \un(\thn)-\ul(\thl) \|^2_{0,\omgb}.
\end{equation}
To show that \eqref{reduced-minimization}  is well-posed we start as in \cite{Gervasio_01_NM,Bochev_14_SINUM} and 
split, for any given $(\thn,\thl)$, the solutions of the state equations into a harmonic and a  homogeneous part. 
The harmonic components  $\vn(\thn)$ and $\vl(\thl)$ of the states solve the equations
\begin{equation}\label{eq:harm}
\left\{\begin{array}{rll}
-\mcL\vn & = 0 & \xb\in\omgn\\ 
\vn & = \thn & \xb\in\omgc\\ 
\vn & = 0 & \xb\in\omgi,
\end{array}\right.
\quad\mbox{and}\quad
\left\{\begin{array}{rll}
-\Delta\vl & = 0 & \xb\in\omgl\\ 
\vl & =\thl & \xb\in\gamc\\ 
\vl & =0 & \xb\in\gami,
\end{array}\right.
\end{equation}
respectively.  The homogeneous components $\un^0$ and $\ul^0$ solve a similar set of equations but with homogeneous volume constraint and boundary condition, respectively:
\begin{equation}\label{eq:homo}
\left\{\begin{array}{rll}
-\mcL\un^0 &= f_n & \xb\in\omgn\\ 
\un^0  &=0 & \xb\in\omgw
\end{array}\right.
\quad\mbox{and}\quad
\left\{\begin{array}{rll}
-\Delta\ul^0 &= f_l & \xb\in\omgl\\ 
\ul^0 & =0 & \xb\in\Gamma_l .
\end{array}\right. 
\end{equation}
In terms of these components
$\un=\vn(\thn)+\un^0$, $ \ul =\vl(\thl)+ \ul^0$, and the objective
\begin{displaymath}
\begin{aligned}
J(\thn,\thl) &
= 
\frac12 \| \vn(\thn)\!-\!\vl(\thl) \|^2_{0,\omgb} + 
\big( \un^0 \!-\! \ul^0, \vn(\thn) \!-\! \vl(\thl) \big)_{0,\omgb}\!\! + 
\frac12 \| \un^0\!-\!\ul^0 \|^2_{0,\omgb}.
\end{aligned}
\end{displaymath}
The Euler-Lagrange equation of  (\ref{reduced-minimization}) is given by:
seek $(\sgn,\sgl)\in \Theta_n\times\Theta_l$ such that
\begin{equation}\label{eq:EL-reduced}
Q(\sgn,\sgl;\mu_n,\mu_l) = F(\mu_n,\mu_l)\quad\forall \, (\mu_n,\mu_l) \in \Theta_n\times\Theta_l,
\end{equation}
where 
$Q(\sgn,\sgl;\mu_n,\mu_l)  = \big(\vn(\sigma_n)-\vl(\sigma_l), \,\vn(\mu_n)-\vl(\mu_l)\big)_{0,\omgb}$ and 
$F(\mu_n,\mu_l)  = -\big( \un^0 - \ul^0, \vn(\mu_n) - \vl(\mu_l) \big)_{0,\omgb}$.
The following lemma establishes a key property of $Q$.

\smallskip\begin{lemma}\label{lemma-ip}
The form $Q(\cdot;\cdot)$ defines an inner product on $\Theta_n\times\Theta_l$.
\end{lemma}
\smallskip{\it Proof.}
By construction $Q(\cdot;\cdot)$ is symmetric and bilinear. Thus, it suffices to show that $Q(\cdot;\cdot)$ is positive definite,
i.e., $Q(\sgn,\sgl; \sgn,\sgl)=0$ if and only if $(\sgn,\sgl)=(0,0)$. Let $(\sgn,\sgl)=(0,0)$ then $\vn(\sgn)=0$ and $\vl(\sgl)=0$, implying $Q(\sgn,\sgl;\sgn,\sgl)=0$. Conversely, if $Q(\sgn,\sgl;\sgn,\sgl)=0$, then $\vn(\sgn)-\vl(\sgl)=0$ in $\omgb$. Let $v = \vn(\sgn)=\vl(\sgl)$ in $\omgb$. 
Then we have that (i) $\Delta v=0$ for all $\xb\in\omgb$, i.e., $v$ is harmonic in $\omgb$, and (ii) $v=0$ for all $\xb\in \omgb\cap\omgi$, i.e., $v$ vanishes on a non-empty interior set of $\omgb$.  By the identity principle for harmonic functions, $v\equiv 0$ in $\overline\omgb$.
Because $\sgn=v$ in $\omgc$ and $\sgl=v$ on $\gamc$ it follows that $(\sgn,\sgl)=(0,0)$.
$\square$

\smallskip As a results, $Q(\cdot;\cdot)$ endows the control space $\Theta_n\times\Theta_l$ with the ``energy'' norm 
\begin{equation}\label{Qnorm} 
\|(\sgn,\sgl)\|^2_* := Q(\sgn,\sgl;\sgn,\sgl) =\|\vn(\sgn)-\vl(\sgl) \|^2_{0,\omgb}.
\end{equation}
Note that $Q$ and $F$ are continuous with respect to the energy norm.
However, the control space $\Theta_n\times\Theta_l$ may not be complete with respect to the energy norm. In this case, following \cite{Abdulle_15,Discacciati_13_SIOPT}, we consider the optimization problem \eqref{reduced-minimization} on the completion $\widetilde\Theta_n\times\widetilde\Theta_l$ of the control space.
\begin{theorem}\label{thm:well-posedness}
Let $\widetilde\Theta_n\times\widetilde\Theta_l$ denote a completion\footnote{If $\Theta_n\times\Theta_l$ is complete, then of course we have that $\widetilde\Theta_n\times\widetilde\Theta_l=\Theta_n\times\Theta_l$.} of the control space with respect to the energy norm \eqref{Qnorm}. Then, the minimization problem
\begin{equation}\label{eq:complete-reduced-minimization}
\min_{(\thn,\thl)\in\widetilde\Theta_n\times\widetilde\Theta_l} J(\thn,\thl)
\end{equation}
has a unique minimizer $(\widetilde\theta_n^*,\widetilde\theta_l^*)\in\widetilde\Theta_n\times\widetilde\Theta_l$ such that
\begin{equation}\label{eq:complete-Euler-Lagrange}
\widetilde Q(\widetilde\theta_n^*,\widetilde\theta_l^*;\mu_n,\mu_l) = \widetilde F(\mu_n,\mu_l)
\quad \forall \, (\mu_n,\mu_l)\in \widetilde\Theta_n\times\widetilde\Theta_l.
\end{equation}
\end{theorem}
{\it Proof.}
Equation \eqref{eq:complete-Euler-Lagrange} is a necessary condition for any minimizer of \eqref{eq:complete-reduced-minimization}. Assume first that the control space is complete, i.e. $\Theta_n\times\Theta_l = \widetilde\Theta_n\times\widetilde\Theta_l$. Then $\Theta_n\times\Theta_l$ is Hilbert and the projection theorem implies that \eqref{eq:complete-Euler-Lagrange} has a unique solution. 

When $\Theta_n\times\Theta_l$ is not complete, the continuous bilinear form $Q$ and the continuous functional $F$ defined on $\Theta_n\times\Theta_l$ can be uniquely extended by using the Hahn-Banach theorem to a continuous bilinear form $\widetilde Q$ and a continuous functional $\widetilde F$ in $\widetilde\Theta_n\times\widetilde\Theta_l$. Then, the existence and uniqueness of the minimizer follow as before. $\square$

To avoid technical distractions, in what follows we assume that the minimizer $(\thn^*,\thl^*)$ belongs to $\Theta_n\times\Theta_l$ and hence 
$\un^*=\un(\thn^*)\in V_{\eta_D}(\omgnp)$ and 
$\ul^*=\ul(\thl^*)\in H^1_\gami(\omgl)$.
We note that in the finite dimensional case the completeness is not an issue, as the discrete control space is Hilbert with respect to the discrete energy norm, see Section \ref{finite-dim}.

\section{Analysis of the LtN coupling error}\label{error-analysis}
We define the LtN coupling error as the $L^2$-norm of the difference between the global nonlocal solution $\wun$ of \eqref{nonlocal-problem} with homogeneous volume constraints and the LtN solution $\unl \in L^2(\omgp)$ given by \eqref{unl-definition}. 
This section shows that the coupling error is bounded by the \emph{modeling error} on the local subdomain, i.e., the error made by replacing the ``true'' nonlocal diffusion operator on $\omgl$ by the Laplacian. 

We prove this result under the following assumptions.
\smallskip

\indent H.1 The kernel $\gam$ satisfies \eqref{eq:gamma-conds} and \eqref{eq:integrable-kernel}.

\indent H.2 The global nonlocal solution\footnote{This assumption can be relaxed to: $\wun\in L^2(\omgp)$ has a well-defined trace on $\gamc$.} $\wun\in H^1(\omgp)$.

\smallskip
We also need the {\it trace} operator $T :H^1({\omgp})\to\Theta_n\!\times\!\Theta_l$ such that
\begin{equation}\label{operatorR}
T(v) := (T_n(v),\;T_l(v)) = (v|_\omgc, \; v|_\gamc)
\qquad\forall \, v\in H^1({\omgp}),
\end{equation}
and the lifting operator  $L:\Theta_n \times \Theta_l\to L^2(\omgp)$, $L(\sgn,\sgl)=H(\sgn,\sgl)+u^0$, where
\begin{equation}\label{operatorB}
H(\sgn,\sgl) := \left\{\begin{array}{ll}
\vn(\sgn) & \omgnp\\ 
\vl(\sgl) & \omgl\setminus\omgb,
\end{array}\right.
\qquad \hbox{and} \qquad
u^0 := \left\{\begin{array}{ll}
\un^0 & \omgnp\\ 
\ul^0 & \omgl\setminus\omgb
\end{array}\right. 
\end{equation}
are a harmonic lifting operator and the homogenous part of the states, respectively. 
Our main result is the following theorem.
\begin{theorem}\label{thm-coupling-error}
Assume that H.1 and H.2 hold. Then, there exists a positive constant $C$ such that
\begin{displaymath}
\|\wun-\unl\|_{0,{\omgp}}\leq C \|\wun-\wul\|_{0,\omgl} \,,
\end{displaymath}
where $\wul = v_l(T_{l}(\wun))+\ul^0$.
\end{theorem}

\smallskip For clarity we break the proof of this theorem into several steps arranged in Sections \ref{sec:prelim}--\ref{sec:proof}.  Also for clarity, below we list  the ancillary results required for the proof of Theorem \ref{thm-coupling-error}, but to avoid clutter, we collect all their proofs in Appendix \ref{app:op-norm}. In what follows we let $C$ and $C_i$, $i=1,2,\ldots$, denote generic positive constants.

\begin{wrapfigure}{r}{0.4\linewidth}
\vspace{-5.4ex}
\centering
\includegraphics[width=0.7\linewidth]{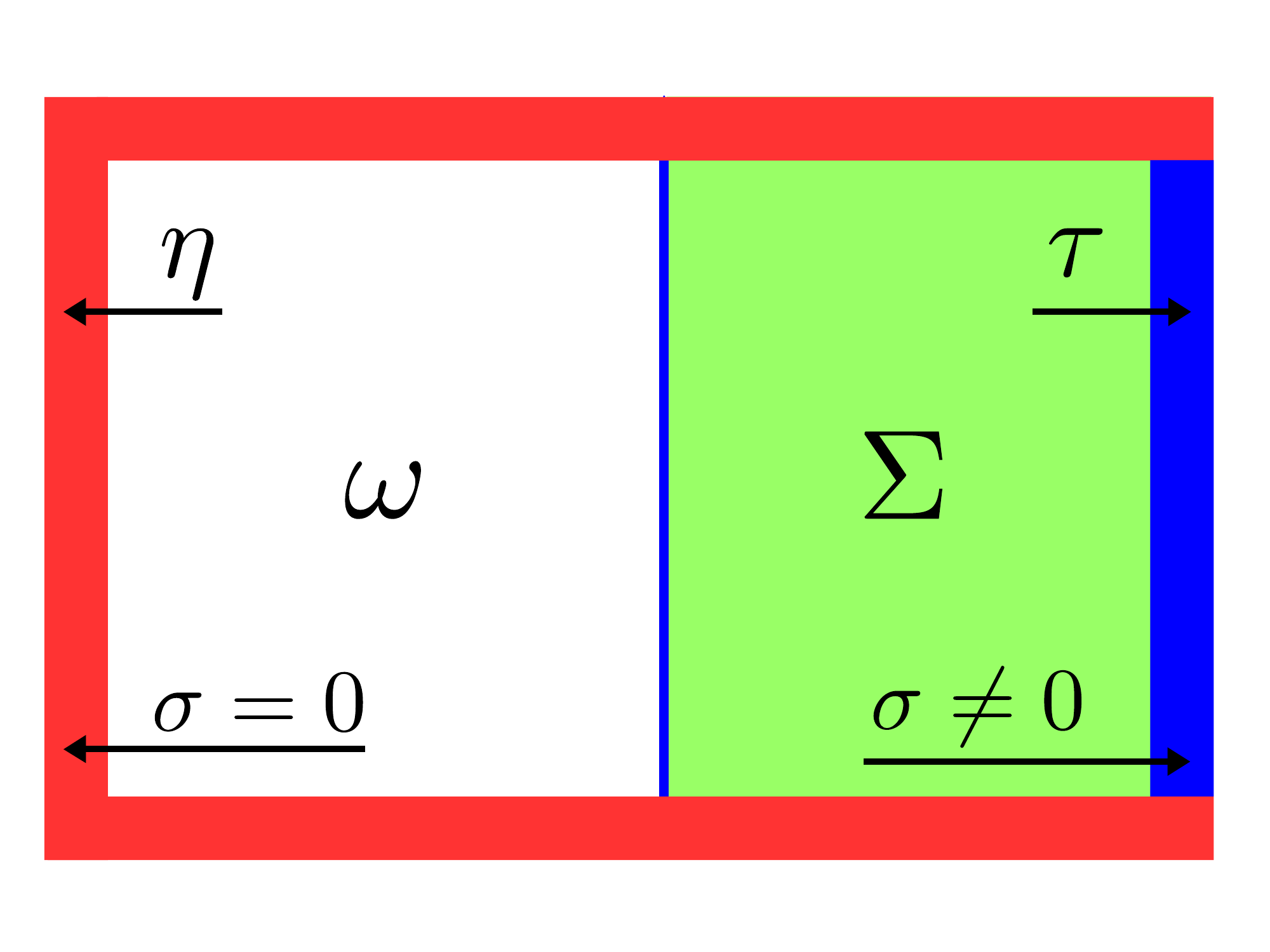}
\noindent
\vspace{-2.8ex}
\caption{Domain for Lemma \ref{lemma-nonlocal-trace}.}
\label{fig:nonlocal-trace}
\end{wrapfigure}

\smallskip
\begin{lemma}\label{lemma-nonlocal-trace}
\emph{[Nonlocal trace inequality]} Let $\omgp=\omega\cup\eta$ where $\omg$ and $\womg$ are an open bounded domain and its associated interaction domain and let $\Sigma\subset{\omgp}$, $\tau\subset\womg$, $\tau\subset\Sigma$ as in Fig. \ref{fig:nonlocal-trace}. Assume that $\sigma\in\widetilde V(\womg)$ is such that $\sigma|_{\womg\setminus\tau}=0$. Then
\begin{equation}\label{nonlocal-trace}
\|\sigma\|_{\widetilde V(\womg)} \leq C |||v |||_{\Sigma}
\quad \forall \, v\in V({\omgp}), \, v|_\womg=\sigma.
\end{equation}
\end{lemma}
\smallskip
\begin{lemma}\label{lem:closed-subsp}
Let $\omgp$, $\omega$ and $\womg$ be defined as in Lemma \ref{lemma-nonlocal-trace}. The trace space $\widetilde V_{\womg\setminus\tau}(\womg)$ is a closed subspace of $L^2(\womg)$.
Furthermore, for all $\mu\in \widetilde V_{\womg\setminus\tau}(\womg)$ we have that 
\begin{equation}\label{eq:trace-space-equivalence}
\|\mu\|_{\widetilde V(\womg)}\leq C\|\mu\|_{0,\womg}.
\end{equation}
\end{lemma}

\smallskip
Application of Lemma \ref{lem:closed-subsp} with $\womg=\omgw$ and $\tau = \omgc$ implies that the trace space $\widetilde V_{\omgi}(\omgw)$ is a closed subspace of $L^2(\omgw)$, and thus it is a Hilbert space in the $L^2$ topology. 
We use this result to prove a strong Cauchy-Schwartz inequality for the nonlocal and local harmonic components of the states, i.e., the solutions to \eqref{eq:harm}. This inequality is essential for the estimate of $\| H\|_{**}$ below.
\smallskip
\begin{lemma}\label{lemma-strongCS}
\emph{[Strong Cauchy-Schwartz inequality]} There exists $\delta<1$ such that for all $(\sgn,\sgl)\in\Theta_n\times\Theta_l$
\begin{equation}\label{c-s}
|(\vn(\sgn),\vl(\sgl))_{0,\omgb}|< \delta \|\vn(\sgn)\|_{0,\omgb} \, \|\vl(\sgl)\|_{0,\omgb}.
\end{equation}
\end{lemma}

\subsection{The harmonic lifting operator is bounded from above}\label{sec:prelim}
We prove that $H$ is bounded in the operator norm $\| \cdot \|_{**}$ induced by the energy norm \eqref{Qnorm}. We refer to Appendix \ref{app:op-norm} for additional notation and auxiliary results used in the proof. We introduce the space 
$\widetilde V_\tau(\eta) = \{\mu\in\widetilde V(\eta): \mu|_\tau=0, \;{\rm for}\;\tau\subseteq\eta\}$.

\smallskip
\begin{lemma}\label{lemma-norm-bound}
Assume that H.1 holds. There exists a positive constant $C<\infty$ such that
\begin{equation}\label{norm-bound}
\|H\|_{**} = \sup\limits_{(0,0)\neq(\sgn,\sgl)\in\Theta_n\times\Theta_l}
\frac{\|H(\sgn,\sgl)\|_{0,\omgp}}{\|(\sgn,\sgl)\|_*}<\dfrac{C}{\kappa},
\end{equation}
where $\kappa$ is the thickness of $\omgb$.
\end{lemma}
{\it Proof.}
To prove \eqref{norm-bound} it suffices to show that
\begin{displaymath}
\|H(\sgn,\sgl)\|_{0,{\omgp}} \leq \widetilde C \|(\sgn,\sgl)\|_* 
\qquad \forall\,(0,0)\neq(\sgn,\sgl)\in\Theta_n\times\Theta_l,
\end{displaymath}
for some positive constant $\widetilde C$, inversely proportional to $\kappa$.
According to the definitions of $H$ and $\|(\cdot,\cdot)\|_*$ in \eqref{operatorB} and \eqref{Qnorm}, this is equivalent to
\begin{displaymath}
\|\chi(\omgnp)v_n(\sgn) + \chi(\omgl\setminus\omgb) v_l(\sgl)\|_{0,{\omgp}}\leq 
\widetilde C\|\vn(\sgn)-\vl(\sgl)\|_{0,\omgb} ,
\end{displaymath}
where $\chi(\cdot)$ is the indicator function. Since $(\omgl\hspace{-.05cm}\setminus\hspace{-.05cm}\omgb)\cap\omgnp=\emptyset$, this inequality reduces to
\begin{displaymath}
\|\vn(\sgn)\|^2_{0,\omgnp} +\|\vl(\sgl)\|^2_{0,\omgl\setminus\omgb}\leq 
\widetilde C^2 \|\vn(\sgn)-\vl(\sgl)\|^2_{0,\omgb}.
\end{displaymath}
The strong Cauchy-Schwarz inequality for the harmonic component (see Lemma \ref{lemma-strongCS}) yields the following lower bound for the right hand side:
\begin{equation}\label{eq:scbs}
\begin{aligned}
      \|\vn(\sgn)- & \vl(\sgl)\|^2_{0,\omgb} = 
      \|\vn(\sgn)\|^2_{0,\omgb} - 2 (\vn(\sgn),\vl(\sgl))_{0,\omgb} 
 +    \|\vl(\sgl)\|^2_{0,\omgb}\\
&\geq \|\vn(\sgn)\|^2_{0,\omgb} - 2 \delta \|\vn(\sgn)\|_{0,\omgb} \, \|\vl(\sgl)\|_{0,\omgb}
 +    \|\vl(\sgl)\|^2_{0,\omgb}\\
&\geq (1-\delta) \big(\|\vn(\sgn)\|^2_{0,\omgb} + \|\vl(\sgl)\|^2_{0,\omgb} \big).
\end{aligned}
\end{equation}
We now proceed to bound $\|\vn(\sgn)\|_{0,\omgb}$ and $\|\vl(\sgl)\|_{0,\omgb}$ from below by $\|\vn(\sgn)\|_{0,\omgnp}$ and $\|\vl(\sgl)\|_{0,\omgl\setminus\omgb}$, respectively. We start with the nonlocal term. 
Let $\mu_n\in\widetilde V_{\omgi}(\omgw)$ denote the extension of $\sgn$ by zero in $\omgi$, i.e., $\mu_n=\sgn$  in $\omgc$ and $0$ in $\omgi$. By using the nonlocal Poincar\'e inequality, the well-posedness of the nonlocal problem and Lemma \ref{lem:closed-subsp} we have
\begin{displaymath}
\begin{array}{rll}
\|\vn(\sgn)\|_{0,\omgnp} \!\! \! & \leq C_{pn} |||\vn(\sgn)|||_\omgnp
& \mbox{nonlocal Poincar\'e inequality \eqref{eq:L2equivalence}} \\ [0.5ex]
&\le C_{pn} K_n \|\mu_n\|_{\widetilde V(\omgw)}
& \mbox{nonlocal well-posedness \eqref{cont-data-nonlocal}} \\ [0.5ex]
& \le C_{pn} K_n K_1 \|\mu_n\|_{0,\eta_n} 
& \mbox{Lemma \ref{lem:closed-subsp}} \\ [0.5ex]
&\le C_{pn} K_n K_1 K_2  \|\vn(\sgn)\|_{0,\omgb}. \!\! & 
\end{array}
\end{displaymath}
Therefore, we have that 
\begin{equation}\label{eq:lbn}
\|\vn(\sgn)\|^2_{0,\omgnp} \leq  \widetilde C_n \|\vn(\sgn)\|^2_{0,\omgb},
\end{equation}
with $\widetilde C_n=C_{pn} K_n K_1 K_2$.

To analyze the local term we derive a Caccioppoli-type inequality for the local harmonic component. We introduce the cutoff function $g\in C^1(\omgl)$ such that $g=1$ in $\overline{\omgl\setminus\omgb}$, $g=0$ in $\Omega\setminus\omgl$, $\|\nabla g\|_{\infty}\leq\frac{1}{\kappa}$, and ${\rm supp}(\nabla g)\subset\omgb$, where $\kappa$ is the thickness of the overlap, see Fig. \ref{domains}. These properties imply that $g\vl$ and $g^2\vl$ belong to $H^1_0(\omgl)$. Next, we note that $\vl$ is the solution of the weak formulation of \eqref{local-problem} with $f_l=0$. Using $g^2\vl$ as a test function then yields the following identity
\begin{displaymath}
0 =  \int_\omgl \nabla\vl \nabla(g^2\vl)\,d\xb 
  = 2\int_\omgl g\vl\,\nabla\vl \nabla g\,d\xb
  +  \int_\omgl g^2\nabla\vl \nabla\vl\,d\xb.
\end{displaymath}
We use the latter to find a bound on $\|\nabla(g\vl)\|_{0,\omgl}$:
\begin{displaymath}
\begin{aligned}
\|\nabla(g\vl)\|^2_{0,\omgl} 
& = \int_\omgl \nabla(g\vl)\nabla(g\vl)\,d\xb - \int_\omgl \nabla\vl \nabla(g^2\vl)\,d\xb \\[2mm]
& = \int_\omgl \nabla(g\vl)\nabla(g\vl)\,d\xb 
  - 2\int_\omgl g\vl\,\nabla\vl \nabla g\,d\xb
  - \int_\omgl g^2\nabla\vl \nabla\vl\,d\xb \\[2mm]
& = \int_\omgl \vl^2 \nabla g \nabla g \,d\xb 
  = \int_\omgb \vl^2 \nabla g \nabla g \,d\xb 
  \leq \frac{1}{\kappa^2} \int_\omgb \vl^2\,d\xb = \frac{1}{\kappa^2} \|\vl\|^2_{0,\omgb}.
\end{aligned}
\end{displaymath}
Thus, we conclude that
\begin{displaymath}
\|\vl\|_{0,\omgl\setminus\omgb} = \|g\vl\|_{0,\omgl\setminus\omgb} \leq \|g\vl\|_{0,\omgl}
\leq C_p \|\nabla(g\vl)\|_{0,\omgl} \leq \frac{C_p}{\kappa} \|\vl\|_{0,\omgb},
\end{displaymath}
where $C_p$ is the local Poincar\'e constant. Let $\widetilde K_{nl}=\max\{\widetilde C_n,C_p^2\}$.
Together with \eqref{eq:scbs} and \eqref{eq:lbn} this yields 
\begin{displaymath}
\begin{aligned}
\|\vn(\sgn)\|^2_{0,\omgnp} +\|\vl(\sgl)\|^2_{0,\omgl\setminus\omgb}
& \leq \dfrac{\widetilde K_{nl}^2}{\kappa^2}\left( \|\vn(\sgn)\|^2_{0,\omgb} +  
       \|\vl(\sgl)\|^2_{0,\omgb} \right) \\
& \leq  \dfrac{\widetilde K_{nl}^2}{\kappa^2(1-\delta)}\|\vn(\sgn)-\vl(\sgl)\|^2_{0,\omgb}.
\end{aligned}
\end{displaymath}
$\square$

\subsection{The approximation error is bounded by the modeling error}\label{sec:approxer}
The optimal solution $(\thna,\thla)$ of the reduced space problem \eqref{reduced-minimization} approximates the trace of the global nonlocal solution $\wun$ on $\omgc$ and $\gamc$, respectively. The following lemma shows that the error in $(\thna,\thla)$ is bounded by the modeling error on $\omgl$.
\smallskip

\begin{lemma}
Let $\wun$ and  $(\thna,\thla)$ solve \eqref{nonlocal-problem} and \eqref{reduced-minimization}, respectively.
Then,
\begin{equation}\label{approx-error}
\|T(\wun)-(\thna,\thla)\|_* \leq \|\wun-\wul\|_{0,\omgb}.
\end{equation}
\end{lemma}
{\it Proof.}
Because $(\thna,\thla)$ satisfies the Euler-Lagrange equation \eqref{eq:EL-reduced} we have
\begin{equation}\label{eq:hEL-reduced}
Q(\thna,\thla;\mu_n,\mu_l) = -\big( \un^0 - \ul^0, \vn(\mu_n) - \vl(\mu_l) \big)_{0,\omgb}\quad\forall \, (\mu_n,\mu_l) \in \Theta_n\times\Theta_l.
\end{equation}
Using this identity together with the energy norm definition \eqref{Qnorm} yields 
\begin{displaymath}
\begin{array}{l}
\displaystyle
Q(T(\wun)-(\thna,\thla); \mu_n,\mu_l)
=
Q(T(\wun); \mu_n,\mu_l)-Q(\thna,\thla; \mu_n,\mu_l) \\[1ex]
\qquad\displaystyle
=
\big(\vn(T_n(\wun))-\vl(T_l(\wun)), \vn(\mu_n) - \vl(\mu_l) \big)_{0,\omgb}+
\big( \un^0 - \ul^0, \vn(\mu_n) - \vl(\mu_l) \big)_{0,\omgb}
\\[1ex]
\qquad\displaystyle
=\big( \wun - \wul, \vn(\mu_n) - \vl(\mu_l) \big)_{0,\omgb}
\le \| \wun - \wul \|_{0,\omgb} \, \| \vn(\mu_n) - \vl(\mu_l) \|_{0,\omgb}
\\[1ex]
\qquad\displaystyle
= \| \wun - \wul \|_{0,\omgb} \, \|(\mu_n, \mu_l)\|_*.
\end{array}
\end{displaymath}
The lemma follows by setting $(\mu_n, \mu_l) = T(\wun)-(\thna,\thla)$ above and observing that
$ Q(T(\wun)-(\thna,\thla); T(\wun)-(\thna,\thla)) = \| T(\wun)-(\thna,\thla) \|^2_* $.
$\square$

\subsection{Proof of Theorem \ref{thm-coupling-error}}\label{sec:proof}
Let $\widehat{u}:=L(T(\wun))$. Definitions \eqref{operatorR} and \eqref{operatorB} together with the identities 
\begin{equation}\label{local_lifting}
\wun|_\omgnp = \vn(T_{n}(\wun))+\un^0
\quad\mbox{and}\quad
\wul = v_l(T_{l}(\wun))+\ul^0,
\end{equation}
imply that
\begin{equation}\label{PRwun}
\widehat{u}=
\left\{\begin{array}{ll}
\wun & \xb\in\omgnp\\ 
\wul & \xb\in\omgl\setminus\omgb.
\end{array}\right.
\end{equation}
Likewise, the identities $\un^*=\vn(\thna) + \un^0$  and $\ul^*=\vl(\thla) + \ul^0$ imply that $\unl = L(\thna,\thla)$.
Adding and subtracting $\widehat{u}$ to the LtN error then yields
\begin{equation}\label{error-splitting}
\begin{array}{rl}
\|\wun-\unl\|_{0,\omgp} & \le \|\wun-\widehat{u}\|_{0,\omgp} + \|\widehat{u} - \unl \|_{0,\omgp} 
\\[0.5ex]
& =  \|\wun-L(T(\wun))\|_{0,\omgp} +  \|L(T(\wun)) - L(\thna,\thla) \|_{0,\omgp}
\\[1ex]
&\displaystyle
  = \|\wun-L(T(\wun))\|_{0,\omgp} 
 + \|H(T(\wun)) - H(\thna,\thla)\|_{0,\omgp} \,.
\end{array}
\end{equation}
The first term in \eqref{error-splitting} is the {\it consistency error} of $L$; \eqref{PRwun} implies that 
\begin{equation}\label{consistency-error}
\|\wun - L(T(\wun))\|_{0,\omgp} = \|\wun - \wul \|_{0,\omgl\setminus\omgb}.
\end{equation}
To estimate the second term we use \eqref{norm-bound} and \eqref{approx-error}:
\begin{displaymath}
\|H(T(\wun)) - H(\thna,\thla)\|_{0,\omgp} \le \| H \|_{**} \|T(\wun)-(\thna,\thla)\|_*
\leq \frac{C}{\kappa}\|\wun-\wul\|_{0,\omgb}.
\end{displaymath}
Combining \eqref{error-splitting} with this bound and \eqref{consistency-error} gives
\begin{displaymath}
\|\wun-\unl\|_{0,\omgp} \leq \|\wun-\wul\|_{0,\omgl\setminus\omgb} 
                        +    \dfrac{C}{\kappa} \|\wun-\wul\|_{0,\omgb} 
                        \leq \left(1+\dfrac{C}{\kappa}\right)\|\wun-\wul\|_{0,\omgl} ,
\end{displaymath}
which completes the proof. $\Box$

\subsection{Convergence of the modeling error}\label{sec:modeling-error-convergence}
In this section we show that $\|\wun-\wul\|_{0,\omgl}$ vanishes as $\varepsilon\to 0$.
\begin{lemma}
Let $\wun$ be the solution of \eqref{nonlocal-problem} with homogeneous volume constraints and let $\wul$ be defined as in \eqref{local_lifting}. Assume that H.1 and H.2 hold; then,
\begin{displaymath}
\lim\limits_{\varepsilon\to 0}\|\wun-\wul\|_{0,\omgl} = 0.
\end{displaymath}
\end{lemma} 
\smallskip
{\it Proof.}
By definition  $\wul$ solves the boundary value problem
\begin{displaymath}
\left\{\begin{array}{rlll}
-\Delta \wul & = & f_l & \xb\in \omgl \\
\wul & = & \wun & \xb\in\gamc \\
\wul & = & 0 & \xb\in\gami,
\end{array}\right.
\end{displaymath}
and so, it is also a solution of the weak equation
\begin{displaymath}
\int_\omgl \nabla\wul\cdot\nabla w d\xb  = \int_\omgl f_l\,w d\xb
\quad\forall \, w\in H^1_0(\omgl).
\end{displaymath}
Let $\psi\in H^1_0(\omgl)$ solve the dual problem
\begin{equation}\label{eq:dual}
\int_{\omgl} \nabla w\cdot\nabla\psi d\xb = \int_{\omgl} (\wul - \wun) w d\xb
\quad\forall\, w\in H^1_0(\omgl)\,.
\end{equation}
Since $\wul - \wun = 0$ on $\Gamma_l$ one can set $w= \wul - \wun$ in \eqref{eq:dual} to obtain 
\begin{displaymath}
\begin{array}{rll}
   \displaystyle \| \wul-\wun \|^2_{0,\omgl} 
=& \displaystyle \int_{\omgl} \nabla (\wul - \wun)\cdot\nabla\psi d\xb  \\[2ex]
=& \displaystyle   \int_{\omgl} f_l \psi d\xb - \int_{\omgl} \nabla \wun\cdot\nabla\psi d\xb 
   \\[2ex]
=& \cancel{\displaystyle  \int_{\omgl\cap\eta} f_l \psi d\xb}
+  \int_{\omgl\setminus\womg} f_n \psi d\xb 
-  \int_{\omgl} \nabla \wun\cdot\nabla\psi d\xb \\[2ex]
=& \displaystyle\int_{\omgl\setminus\womg} \mcL\wun \psi d\xb 
-  \int_{\omgl\setminus\womg} \nabla \wun\cdot\nabla\psi d\xb\to 0, \;\;\quad
   \mbox{as $\varepsilon\to 0$,}
\end{array}
\end{displaymath}
where the third equality follows from the fact that $f_l$ is extended to zero in $\eta$ and the limit follows from the result in \cite[Section~5]{Du_13_MMMAS}\footnote{It can be shown that when kernels satisfying \eqref{eq:gamma-conds} and \eqref{eq:integrable-kernel} are properly scaled, so that $\mcL\to\Delta$, $\|\widehat u_n-\widehat u_l\|_{0,\omgl}\leq C\varepsilon^2+\mcO(\varepsilon^4)$. The same result holds for the peridynamics kernel in \eqref{test-kernel}.}. 
$\square$

\begin{remark}
An immediate consequence of the results in this section is the fact that the LtN solution converges to the corresponding local solution in the limit of vanishing nonlocality. This means that, if an asymptotically-compatible (AC) discretization is used, the discreticretized LtN solution converges to the corresponding continuous local solution as $\epsilon\to 0$ and $h\to 0$ simultaneously, where $h$ is the discretization parameter. In other words, the proposed formulation yields an AC scheme, provided AC discretizations are employed.
\end{remark}

\section{Approximation of the optimization-based LtN formulation} \label{finite-dim}
This section presents the discretization and the error analysis of the LtN formulation \eqref{minimization}. Throughout this section we assume that $\omgc$ and $\gamc$ are polygonal domains; this assumption is not restrictive as those are virtual domains that we can define at our discretion.
\subsection{Discretization}\label{sec:discrete}
We use a reduced-space approach to solve the optimization-based LtN problem \eqref{minimization}  numerically, i.e., we discretize and solve the problem
\begin{equation}\label{eq:rsp}
\displaystyle\min\limits_{\thn,\thl} J(\thn,\thl) 
\quad {\rm with} \quad
J(\thn,\thl)  = \frac12 \| \un(\thn)-\ul(\thl) \|^2_{0,\omgb} 
\end{equation}
where $\un(\thn)\in V_\omgi(\omgnp) $ solves the weak nonlocal equation
\beq\label{var-state-nn}
      {B_n} (\un(\thn); \zn,\kappa_n ) 
:=\!  \int\limits_\omgnp \! \int\limits_\omgnp \mcG\un\cdot\mcG\zn\,d\yb d\xb
+ \! \int\limits_\omgc \un \kappa_n  \,d\xb 
= \! \int\limits_\omgn \! f_n\zn\,d\xb +\!\int\limits_\omgc  \! \thn \kappa_n \,d\xb,
\eeq
for all $(\zn,\kappa_n)\in V_\omgw\times\Theta_n^*$, and $\ul(\thl)\in H^1_{\gami}(\omgl)$ solves the weak local equation 
\beq\label{var-state-ll}
   {B_l}(\ul(\thl); \zl,\kappa_l ) 
:= \int_\omgl \nabla\ul\nabla\zl\,d\xb + \int_\gamc  \ul \kappa_l \,d\xb 
 = \int_\omgl f_l\zl\,d\xb + \int_\gamc \thl \kappa_l \,d\xb,
\eeq
for all $(\zl,\kappa_l)\in H^1_0(\omgl)\times \Theta_l^*$. Here, $\Theta_n^*$ and $\Theta_l^*$ are the duals of $\Theta_n$ and $\Theta_l$, respectively.  

To discretize \eqref{eq:rsp}--\eqref{var-state-ll} we consider the following conforming finite element spaces \cite{DElia_2020Acta,DElia_2020Cookbook}
\begin{equation}\label{finite-dim-spaces}
\begin{array}{lll}
V_\omgi^h\subset V_\omgi(\omgnp), &\;\; 
\Theta^h_n \subset \Theta_n, &\,\,
V_\omgw^h\times\Theta^h_n \,\subset \, V_\omgw(\omgnp) \times \Theta_n^*, \\[3mm]
H^h_\gami\subset H^1_\gami(\omgl), &\;\; 
\Theta^h_l \subset \Theta_l, &\,\,
H^h_0\times\Theta^h_l \,\,\subset \,\, H^1_0(\omgl) \times \Theta_l^*
\end{array}
\end{equation}
for the nonlocal and local states, controls, and test functions\footnote{For simplicity, we approximate $\Theta_n$, $\Theta_l$ and their duals by the same finite dimensional space.}, respectively. In general, the finite element spaces for the nonlocal and local problems can be defined on different meshes with parameters $h_n>0$ and $h_l>0$, respectively, and can have different polynomial orders given by integers  $p_n\ge 1$ and $p_l\ge 1$, respectively. 

Restriction of  \eqref{eq:rsp}--\eqref{var-state-ll}  to the finite element spaces \eqref{finite-dim-spaces} defines the discrete reduced-space LtN formulation
\begin{equation}\label{eq:rsph}
\displaystyle\min\limits_{\thnh,\thlh} J_h(\thnh,\thlh)
\quad {\rm with} \quad
J_h(\thnh,\thlh)=\frac12 \| \unh(\thnh)-\ulh(\thlh) \|^2_{0,\omgb} 
\end{equation}
where $\unh(\thnh)\in V^h_\omgi$ solves the discrete nonlocal state equation
\beq\label{var-state-nnh}
{B_n}(\unh(\thnh); \znh,\kappa_n^h) = \int_\omgn f_n\znh\,d\xb +  \int_\omgc \thnh\kappa_n^h\,d\xb,
\quad\mbox{$\forall \; (\znh,\kappa_n^h)\in V^h_\omgw\times \Theta_n^h$},
\eeq
and $\ulh(\thlh)\in H^h_\gami$ solves the discrete local state equation\footnote{Note that both \eqref{var-state-nnh} and \eqref{var-state-llh} are well-posed.}
\beq\label{var-state-llh}
{B_l}(\ulh(\thlh); \zlh,\kappa_l^h)=\int_\omgl f_l\zl\,d\xb + \int_\gamc \thlh\kappa_l^h\,d\xb,
\quad\mbox{$\forall\;(\zlh,\kappa_l^h)\in H^h_0\times \Theta^h_l$}.
\eeq 
Following Section \ref{well-posedness}, we write the solutions of \eqref{var-state-nnh} and \eqref{var-state-llh} as
\begin{equation}\label{eq:discr-state-decomp}
\unh = \vnh+\unoh \quad {\rm and} \quad \ulh = \vlh + \uloh,
\end{equation}
where $\vnh$ and $\vlh$ are ``harmonic'' components solving \eqref{var-state-nnh} and \eqref{var-state-llh} with $f_n=0$ and $f_l=0$ respectively, whereas $\unoh$ and $\uloh$ are ``homogeneous'' components  solving \eqref{var-state-nnh} and \eqref{var-state-llh} with $\thnh=0$ and $\thlh=0$, respectively. In terms of these components
\begin{displaymath}
J_h(\thnh,\thlh) = 
\frac12 \| \vnh(\thnh)\!-\!\vlh(\thlh) \|^2_{0,\omgb} + 
\big( \unoh \!-\! \uloh, \vnh(\thnh) \!-\! \vlh(\thlh) \big)_{0,\omgb}\!\! + 
\frac12 \| \unoh \!-\!\uloh \|^2_{0,\omgb} \,.
\end{displaymath}
The Euler-Lagrange equation of \eqref{eq:rsph} has the form:
seek $(\sgnh,\sglh)\in \Theta^h_n\times\Theta^h_l$ such that
\begin{equation}\label{eq:disEL}
Q_h(\sgnh,\sglh;\mu^h_n,\mu^h_l) = F_h(\mu^h_n,\mu^h_l)\quad\forall \,(\mu^h_n,\mu^h_l) 
\in \Theta^h_n\times\Theta^h_l,
\end{equation}
\smallskip where 
$Q_h(\sgnh,\sglh;\mu^h_n,\mu^h_l)  = \big(\vnh(\sgnh)-\vlh(\sglh), 
\,\vnh(\mu^h_n)-\vlh(\mu^h_l)\big)_{0,\omgb}$
and
$
F_h(\mu^h_n,\mu^h_l) = -\big( \unoh - \uloh, \vnh(\mu^h_n) - \vlh(\mu^h_l) \big)_{0,\omgb}$.
To prove the positivity of $Q_h$, the arguments of Lemma \ref{lemma-ip} cannot be extended, as the identity principle does not hold for $\vlh$. We use instead the discrete strong Cauchy-Schwarz inequality in Lemma \ref{lemma-discrete-SCS}.

\smallskip
\begin{lemma}\label{lemma-d-ip}
The form $Q_h(\cdot,\cdot)$ defines an inner product on $\Theta_n^h\times\Theta_l^h$.
\end{lemma}

\smallskip
{\it Proof.}
We prove that  $Q_h(\sgnh,\sglh; \sgnh,\sglh)=0$ if and only if $(\sgnh,\sglh)=(0,0)$.
If $(\sgnh,\sglh)=(0,0)$ then $\vnh(\sgnh)=0$ and $\vlh(\sglh)=0$, implying $Q_h(\sgnh,\sglh; \sgnh,\sglh)=0$. 
Conversely, if $Q_h(\sgnh,\sglh; \sgnh,\sglh)=0$, then 
\begin{displaymath}
0=Q_h(\sgnh,\sglh; \sgnh,\sglh)= \|\vnh(\sgnh)\|_{0,\omgb}^2 + \|\vlh(\sglh)\|^2_{0,\omgb} -2(\vnh(\sgnh),\vlh(\sglh))_{0,\omgb}.
\end{displaymath}
The discrete strong Cauchy-Schwarz inequality (see Lemma \ref{lemma-discrete-SCS}) then implies
\begin{equation}\label{discrete_coer}
0\geq(1-\delta) (\|\vnh(\sgnh)\|^2_{0,\omgb} + \|\vlh(\sglh)\|^2_{0,\omgb}).
\end{equation}
Since $\delta<1$ the left hand side in the above inequality is nonnegative. Thus, we must have that  $\vnh(\sgnh)=0$ and $\vlh(\sglh)=0$, which implies $(\sgnh,\sglh)=(0,0)$.
$\square$

Lemma \ref{lemma:HilbertQh} proves that  $\Theta^h_n\times\Theta^h_l$ is Hilbert with respect to the discrete energy norm
\begin{equation}\label{eq:dis-en}
 \|\mu^h_n,\mu^h_l\|^2_{h*}:=Q_h(\mu^h_n,\mu^h_l;\mu^h_n,\mu^h_l).
 \end{equation}
This fact, Lemma \ref{lemma-d-ip} and the projection theorem provide the following corollary.
\begin{corollary}\label{cor:rsph}
The reduced space problem \eqref{eq:rsph} has a unique minimizer.
\end{corollary}

\subsection{Convergence analysis}\label{sec:conv-anl}
In this section we prove that the discrete solution $(\thn^{h*},\thl^{h*})$ converges to the exact solution $(\thn^*,\thl^*)$ assuming the latter belongs to the ``raw'' control space $\Theta_n\times\Theta_l$. This assumption mirrors the one made in \cite{Abdulle_15} and is necessary because the continuous problem is well-posed in the completion of the raw control space. We prove this result under the following assumptions.
\smallskip

\indent H.3 The optimal solution belongs to the raw space: $(\thn^*,\thl^*)\in \Theta_n\times\Theta_l$.

\indent H.4 The kernel $\gamma$ is translation invariant, i.e. $\gamma(\xb,\yb)=\gamma(\xb-\yb)$\footnote{Note that this assumption is not too restrictive; in fact, it is very common in nonlocal mechanics applications.}.

\smallskip
Let $(\thn^*,\thl^*)\in \Theta_n\times\Theta_l$ denote the optimal solution of \eqref{reduced-minimization} and $(\thn^{h*},\thl^{h*})\in \Theta^h_n\times\Theta^h_l$ be the optimal solution of its discretization \eqref{eq:rsph}. We denote the associated optimal states by $(\un^*,\ul^*)$, and $(\un^{h*},\ul^{h*})$, respectively, that is,
\begin{displaymath}
(\un^*, \ul^*) = (\un(\thn^*), \ul(\thl^*))
\quad\mbox{and}\quad
(\un^{h*}, \ul^{h*}) = (\un^h(\thn^{h*}),\ul(\thl^{h*})).
\end{displaymath}
We will estimate the discrete energy norm of the error $(\thn^*-\thn^{h*};\thl^*-\thl^{h*})$ using Strang's second\footnote{The discrete problem \eqref{eq:disEL} also fits in the setting of Strang's first Lemma \cite[Lemma 2.27, p.95]{Ern_04_BOOK}. We use the second lemma because it simplifies the analysis.}
Lemma; see, e.g., \cite[Lemma 2.25, p.94]{Ern_04_BOOK}. Application of this lemma is contingent upon two conditions: (i) the discrete form $Q_h$ is continuous and coercive with respect to $\|\cdot\|_{h*}$, and (ii) there exists a positive real constant $C$ such that 
\begin{equation}\label{eq:strang-ass}
\| \mu_n,\mu_l \|_{h*} \le C \| \mu_n,\mu_l \|_{*}\quad
\forall \, (\mu_n,\mu_l)\in \Theta_n\times\Theta_l \,.
\end{equation}
The first assumption holds trivially. To verify \eqref{eq:strang-ass} note that 
$$
\| \mu_n,\mu_l \|_{h*}^2  := Q_h(\mu_n,\mu_l;\mu_n,\mu_l) = \| \vnh(\mu_n)-\vlh(\mu_l)\|_{0,\omgb}^2 \,.
$$
Given $\mu_l\in \Theta_l$ the function $\vl^h(\mu_l)$ solves the weak equation 
$$
{B_l}(\vlh(\mu_l); \zlh,\kappa_l^h)=\int_\gamc \mu_l\kappa_l^h\,d\xb
=
\int_\gamc \Pi_l(\mu_l)\kappa_l^h\,d\xb
\quad\mbox{$\forall\;(\zlh,\kappa_l^h)\in H^h_0\times \Theta^h_l$},
$$
where $\Pi_l$ is the $L^2$ projection onto $\Theta^h_l$, i.e.,  $\vl^h(\mu_l) =  \vl^h(\Pi_l\mu_l)$. Similarly, we have that $\vn^h(\mu_n) =  \vn^h(\Pi_n\mu_n)$, where $\Pi_n$ is the $L^2$ projection onto $\Theta^h_n$. 
%
%
Additionally, similarly to \cite{Abdulle_15}, we assume that there exist positive constants $\gamma^*_n$, $\gamma_l^*$, $\gamma_{n*}$, and $\gamma_{l*}$ such that for $h_n$ and $h_l$ small enough the following inequalities hold:
\begin{equation}\label{eq:discr-cont-equivalence}
\begin{aligned}
\gamma_{n*} \|\vn(\sgnh)\|_{0,\omgb} & \leq \|\vnh(\sgnh)\|_{0,\omgb} 
                                     & \leq \gamma_n^*\|\vn(\sgnh)\|_{0,\omgb}\\ 
\gamma_{l*} \|\vl(\sglh)\|_{0,\omgb} & \leq \|\vlh(\sglh)\|_{0,\omgb} 
                                     & \leq  \gamma_l^*\|\vl(\sglh)\|_{0,\omgb}.
\end{aligned}
\end{equation}
The latter, the strong Cauchy-Schwarz inequality and the boundedness of the $L^2$ projection operators yield
\begin{equation}
\begin{aligned}
       \| \mu_n,\mu_l \|_{h*}^2 
& =    \| \vnh(\mu_n)-\vlh(\mu_l)\|_{0,\omgb}^2 
  =    \| \vnh(\Pi_n(\mu_n))-\vlh(\Pi_l(\mu_l))\|_{0,\omgb}^2 \\[2ex] 
& \leq \| \vnh(\Pi_n(\mu_n))\|_{0,\omgb}^2 + \|\vlh(\Pi_l(\mu_l))\|_{0,\omgb}^2 \\[2ex]
& \leq \gamma_n^*\|\vn(\Pi_n(\mu_n))\|_{0,\omgb}^2 + \gamma_l^*
       \|\vl(\Pi_l(\mu_l))\|_{0,\omgb}^2 \\[2ex]
& \leq \frac{C_1}{1-\delta}\| \vn(\Pi_n(\mu_n))-\vl(\Pi_l(\mu_l))\|_{0,\omgb}^2 \\[2ex]
& =    \frac{C_1}{1-\delta}\|\Pi_n\mu_n,\Pi_l\mu_l \|_{*}^2
  \leq \frac{C_2}{1-\delta}\|\mu_n,\mu_l \|_{*}^2.
\end{aligned}
\end{equation}
Application of Strang's second lemma then yields the following error estimate.
\begin{lemma}\label{error}
Let $(\thna,\thla)$ and $(\thn^{h*},\thl^{h*})$ be the solutions of \eqref{eq:EL-reduced} and \eqref{eq:disEL}, then 
\begin{equation}\label{eq:strang-EL}
\begin{aligned}
& \|(\thna-\thn^{h*},  \thl^*-\thl^{h*})\|_{h*} \leq \\ 
& \inf_{(\sgn^h,\sgl^h)} \|(\thn^*-\sgn^h,\thl^*-\sgl^h)\|_{h*} 
+ 
\sup_{\|(\mu^h_n,\mu^h_l)\|_{h*}=1}  |Q_h(\thn^*,\thl^*;\mu^h_n,\mu^h_l) - F_h(\mu^h_n,\mu^h_l)|
\end{aligned}
\end{equation}
where $(\sgn^h,\sgl^h),\;(\mu^h_n,\mu^h_l)\in \Theta_n^h\times\Theta_l^h$.
\end{lemma} 

\smallskip
We use the result in Lemma \ref{error} to obtain asymptotic convergence rates under the assumption that 
1) the homogeneous problems \eqref{eq:homo} have solutions $\un^0\in H^{p_n+t}(\omgnp)$, for $t\in[0,1]$, and $\ul^0\in H^{p_l+1}(\omgl)$; 
2) the control variables $(\thn,\thl)$ are such that $\un(\thn)\in H^{p_n+t}(\omgnp)$ and $\ul(\thl)\in H^{p_l+1}(\omgl)$.
We treat the first term in \eqref{eq:strang-EL} by using the norm-equivalence \eqref{eq:energy-equiv}; we have
\begin{equation}\label{eq:strag-first}
\begin{aligned}
         \inf_{(\sgn^h,\sgl^h)} \|(\thn^*-\sgn^h,&\thl^*-\sgl^h)\|^2_{h*} 
     =   \inf_{(\sgn^h,\sgl^h)} \|(\Pi_n\thn^*-\sgn^h,\Pi_l\thl^*-\sgl^h)\|^2_{h*} \\ 
& \leq C \inf_{(\sgn^h,\sgl^h)}
         \|(\Pi_n\thn^*-\sgn^h,\Pi_l\thl^*-\sgl^h)\|^2_{\Theta_n\times\Theta_l} = 0,
\end{aligned}
\end{equation}
where $\|\cdot\|_{\Theta_n\times\Theta_l}$ is defined as in \eqref{eq:control-norm}. We focus on the second term in \eqref{eq:strang-EL}. Adding and subtracting $Q(\thn^*,\thl^*;\mu^h_n,\mu^h_l)$ and using conformity of $\Theta^h_n\times\Theta^h_l$ gives
\begin{displaymath}
\begin{array}{l}
Q_h(\thn^*,\thl^*;\mu^h_n,\mu^h_l) - F_h(\mu^h_n,\mu^h_l)=
 \\ [1.5ex]
\qquad
=\big[Q_h(\thn^*,\thl^*;\mu^h_n,\mu^h_l) -Q(\thn^*,\thl^*;\mu^h_n,\mu^h_l)\big] 
+
\big[Q(\thn^*,\thl^*;\mu^h_n,\mu^h_l)- F_h(\mu^h_n,\mu^h_l)\big]  \\[1.5ex] 
\qquad
=\big[Q_h(\thn^*,\thl^*;\mu^h_n,\mu^h_l) -Q(\thn^*,\thl^*;\mu^h_n,\mu^h_l)\big] 
+
\big[F(\mu^h_n,\mu^h_l)- F_h(\mu^h_n,\mu^h_l)\big] .
\end{array}
\end{displaymath}
Adding and subtracting the terms
\begin{displaymath}
(\vn(\thn^*)-\vl(\thl^*),\vnh(\mu^h_n)-\vlh(\mu^h_l))_{0,\omgb}
\quad\mbox{and}\quad
(\un^0-\ul^0,\vnh(\mu^h_n)-\vlh(\mu^h_l))_{0,\omgb}
\end{displaymath}
to the last expression and using the definitions of $Q$, $F$, $Q_h$ and $F_h$ yields the identity:
\begin{displaymath}
\begin{array}{l}
Q_h(\thn^*,\thl^*;\mu^h_n,\mu^h_l) - F_h(\mu^h_n,\mu^h_l)= \\[1ex]
\qquad
= \left((\vnh(\thn^*)-\vn(\thn^*))-(\vlh(\thl^*)-\vl(\thl^*)),\vnh(\mu^h_n)-\vlh(\mu^h_l)\right)_{0,\omgb} +\\[1ex]
\qquad 
\left((\un^{h0}-\un^0)-(\ul^{h0}-\ul^0),\vnh(\mu^h_n)-\vlh(\mu^h_l)\right)_{0,\omgb} + \\[1ex]
\qquad 
\left(\un^*-\ul^*, (\vnh(\mu^h_n)-\vn(\mu^h_n) - (\vlh(\mu^h_l) - \vl(\mu^h_l)\right)_{0,\omgb}. 
\end{array}
\end{displaymath}
Application of the Cauchy-Schwartz inequality then gives the following upper bound:
\begin{displaymath}
\begin{array}{l}
\big| Q_h(\thn^*,\thl^*;\mu^h_n,\mu^h_l) - F_h(\mu^h_n,\mu^h_l) \big| \le \| \vnh(\mu^h_n)\!-\!\vlh(\mu^h_l)\|_{0,\omgb} \times\\ [1ex]
\qquad
\Big(
\| \vnh(\thn^*)\!-\!\vn(\thn^*)\|_{0,\omgb} \!+\!
\|\vlh(\thl^*)  \!-\!   \vl(\thl^*)\|_{0,\omgb} \!+\!
\| \un^{h0}\!-\!\un^0 \|_{0,\omgb} \!+\! 
\| \ul^{h0}\!-\!\ul^0 \|_{0,\omgb}
\Big)  +\\[1ex]
\qquad 
\Big( 
\|\vnh(\mu^h_n)-\vn(\mu^h_n) \|_{0,\omgb}  +
\|\vlh(\mu^h_l) - \vl(\mu^h_l) \|_{0,\omgb} 
\Big) \times 
\| \un^{*}\!-\!\ul^{*} \|_{0,\omgb}\,.
\end{array}
\end{displaymath}
Furthermore, note that 
$\| \vnh(\mu^h_n)\!-\!\vlh(\mu^h_l)\|_{0,\omgb} = \|\mu^h_n, \mu^h_l \|_{h*} = 1$,
and that $\| \un^{*}\!-\!\ul^{*} \|_{0,\omgb} = J(\un^*,\ul^*)$ is the optimal value of the objective functional, which is bounded by the modeling error. 
The regularity assumptions on the nonlocal solutions in \eqref{eq:homo} allow us to apply Theorem 6.2 in \cite[p.689]{Du_12_SIREV}: 
\begin{displaymath}
\| \un^{h0}\!-\!\un^0 \|_{0,\omgb}  \le C\,h^{p_n+t}_n \| \un^0 \|_{p_n+t,\omgnp},
\end{displaymath}
where $t\in[0,1]$.
Furthermore, the regularity assumptions on the local solutions in \eqref{eq:homo} allow us to use Corollary 1.122 in \cite[p.66]{Ern_04_BOOK} to conclude that
\begin{displaymath}
\| \ul^{h0}\!-\!\ul^0 \|_{0,\omgb} \le Ch^{p_l+1}_l \| \ul^0 \|_{p_l+1,\omgl}.
\end{displaymath}
According to Weyl's Lemma \cite{Weyl_40_DMJ} the local harmonic liftings $\vl(\thl^*)$ and $\vl(\mu^h_l)$ are smooth functions and so there are positive constants $C_1$ and $C_2$ such that 
\begin{displaymath}
\|\vlh(\thl^*)  \!-\!   \vl(\thl^*)\|_{0,\omgb} \leq C_1 h^{p_l+1}_l
\quad \hbox{and} \quad
\|\vlh(\mu^h_l) - \vl(\mu^h_l) \|_{0,\omgb} \leq C_2 h^{p_l+1}_l.
\end{displaymath}
While a similar result holds for the nonlocal harmonic lifting $\vn(\thn^*)$, the treatment of $\vn(\mu_n^h)$ is more involved, due to the discrete nature of the Dirichlet data, and it requires 
an auxiliary function $\widetilde \mu_n\in C^\infty(\omgc)$ such that
$\|\widetilde\mu_n-\mu_n^h\|_{L^2(\omgc)}\leq \epsilon$, for an arbitrarily small $\epsilon$. 
Because $\vn$ and $\vnh$ depend continuously on the data, 
\begin{equation}\label{eq:nonlocal-harm-lift}
\begin{aligned}
\|&\vnh(\mu^h_n)  - \vn(\mu^h_n) \|_{0,\omgb}                      \\ 
& \leq \|\vnh(\mu^h_n)-\vnh(\widetilde\mu_n) \|_{0,\omgb} 
  +    \|\vnh(\widetilde\mu_n)-\vn(\widetilde\mu_n) \|_{0,\omgb}   
  +    \|\vn(\widetilde\mu_n)-\vn(\mu^h_n) \|_{0,\omgb}             \\ 
& \leq C_1 h_n^{p_n+t}\|\vn(\widetilde\mu_n)\|_{p_n+t,\omgnp}
  +    \|\vnh(\mu^h_n-\widetilde\mu_n) \|_{0,\omgb} 
  +    \|\vn(\widetilde\mu_n-\mu^h_n) \|_{0,\omgb}                  \\ 
& \leq C_1 h_n^{p_n+t}\|\vn(\widetilde\mu_n)\|_{p_n+t,\omgnp}
  +    C_2 \|\mu^h_n-\widetilde\mu_n\|_{0,\omgc} 
  +    C_3 \|\widetilde\mu_n-\mu^h_n\|_{0,\omgc}.  
\end{aligned}
\end{equation}
Since $\epsilon$ can be arbitrarily small, the last two terms in \eqref{eq:nonlocal-harm-lift} are negligible. 
To complete the estimate we only need a uniform bound on $\|\vn(\widetilde\mu_n)\|_{p_n+t,\omgnp}$. To this end, assume that for all $\widetilde\mu_n\in C^\infty(\omgc)$, $\vn(\widetilde\mu_n)\in C^k(\omgnp)$ with $k= p_n+t$. Under this assumption 
${\rm D}^\beta[\mcL \vn(\widetilde\mu_n)]=\mcL {\rm D}^\beta[\vn(\widetilde\mu_n)]$ 
for all $\beta\leq k$. 
Taking into account that $\vn(\widetilde\mu_n)$ is nonlocal harmonic, i.e., $\mcL \vn(\widetilde\mu_n)=0$, it follows that $\mcL {\rm D}^\beta[\vn(\widetilde\mu_n)]=0$, i.e., ${\rm D}^\beta[\vn(\widetilde\mu_n)]$ is also nonlocal harmonic for all $\beta\leq k$. Thus, ${\rm D}^\beta[\vn(\widetilde\mu_n)]$ has a uniformly bounded $L^2$ norm, i.e. $\|{\rm D}^\beta[v_n]\|_{0,\omgnp} \leq C_\beta, \; \forall\, \beta\leq k$. This implies the existence of a positive constant $C$ such that, $\|\vn(\widetilde\mu_n)\|_{p_n+t,\omgnp}\leq C$.
It follows that there exist positive constants $C_1$ and $C_2$ such that
\begin{displaymath}
\|\vnh(\thn^*)  \!-\!   \vn(\thn^*)\|_{0,\omgb}\leq C_1 h^{p_n+t}_n
\quad \hbox{and} \quad
\|\vnh(\mu^h_n) - \vn(\mu^h_n) \|_{0,\omgb} \leq C_2 h^{p_n+t}_n.
\end{displaymath}
We have just shown the following result.

\begin{theorem}\label{thm:hstar-discr-error}
Assume that H.1--H.4 hold. Then, there exist positive constants $C_1,\,C_2$ such that
\begin{equation}\label{eq:hstar-discr-error}
\|(\thna-\thn^{h*},  \thl^*-\thl^{h*})\|_{h*}
\leq C_1 h_n^{p_n+t} + C_2 h_l^{p_l+1}.
\end{equation}
\end{theorem}
We use Theorem \ref{thm:hstar-discr-error} to estimate the $\Theta_n\times\Theta_l$ norm of the discretization error.
\begin{corollary}\label{cor:L2H12-discr-error}
Assume that H.1--H.4 hold. Then, there exist positive constants $C_1,\,C_2,\,C_3$ such that
\begin{equation}\label{eq:L2H12-discr-error}
\|(\thna-\thn^{h*},  \thla-\thl^{h*})\|^2_{\Theta_n\times\Theta_l}  
\leq C_1 h_n^{2(p_n+t)} + C_2 h_l^{2p_l+1}.
\end{equation}
\end{corollary}

{\it Proof.}
Adding and subtracting $\Pi_n\thna$ and $\Pi_l\thla$, and using the triangle inequality
\begin{equation}\label{eq:error-splitting}
\begin{aligned}
       \|(\thna&-\thn^{h*}, \thla-\thl^{h*})\|_{\Theta_n\times\Theta_l} \\ 
& \leq \|(\thna-\Pi_n\thna,  \thla-\Pi_l\thla)\|_{\Theta_n\times\Theta_l} 
     + \|(\Pi_n\thna-\thn^{h*},  \Pi_l\thla-\thl^{h*})\|_{\Theta_n\times\Theta_l}.
\end{aligned}
\end{equation}
Using standard finite element approximation results for the first term yields
\begin{equation}\label{eq:error-splitting-first}
     \|(\thna-\Pi_n\thna,  \thla-\Pi_l\thla)\|^2_{\Theta_n\times\Theta_l} 
\leq C_2 h_n^{2(p_n+t)} \|\thn^*\|^2_{p_n+t,\omgc}
   + C_3 h_l^{2p_l+1} \|\thl^*\|^2_{p_l+\frac12,\gamc}.
\end{equation}
We focus on the second term in \eqref{eq:error-splitting}. By the norm-equivalence \eqref{eq:energy-equiv} in the discrete control space, we have
\begin{displaymath}
\|(\Pi_n\thna-\thn^{h*},  \Pi_n\thna-\thl^{h*})\|_{\Theta_n\!\times\Theta_l} 
\!  \leq  \! C \|(\Pi_n\thna-\thn^{h*},  \Pi_n\thna-\thl^{h*})\|_{h*}
    \!= \!C \|(\thna-\thn^{h*}, \thna-\thl^{h*})\|_{h*}.
\end{displaymath}
This result along with \eqref{eq:hstar-discr-error} and \eqref{eq:error-splitting-first} implies \eqref{eq:L2H12-discr-error}.
$\square$

Since $\un^*$ and $\ul^*$ depend continuously on the data, \eqref{eq:L2H12-discr-error}, Corollary \ref{cor:L2H12-discr-error} implies that
\begin{equation}\label{eq:states-discr-error}
\begin{aligned}
\|\un^*-\un^{*h}\|^2_{0,\omgnp} & \leq K_{n1} \, h_n^{2(p_n+t)}+ 
                                       K_{n2} \, h_l^{2p_l+1} \\ 
\|\ul^*-\ul^{*h}\|^2_{0,\omgl}  & \leq K_{l1} \, h_n^{2(p_n+t)}+ 
                                       K_{l2} \, h_l^{2p_l+1},
\end{aligned}
\end{equation}
that is, the $L^2$ norm error of the state variables is of the same order as the $L^2\times H^\frac12$ norm error of the controls.

\section{Numerical tests}\label{numerical-tests}
\begin{figure}
\centering
\includegraphics[width=0.8\linewidth]{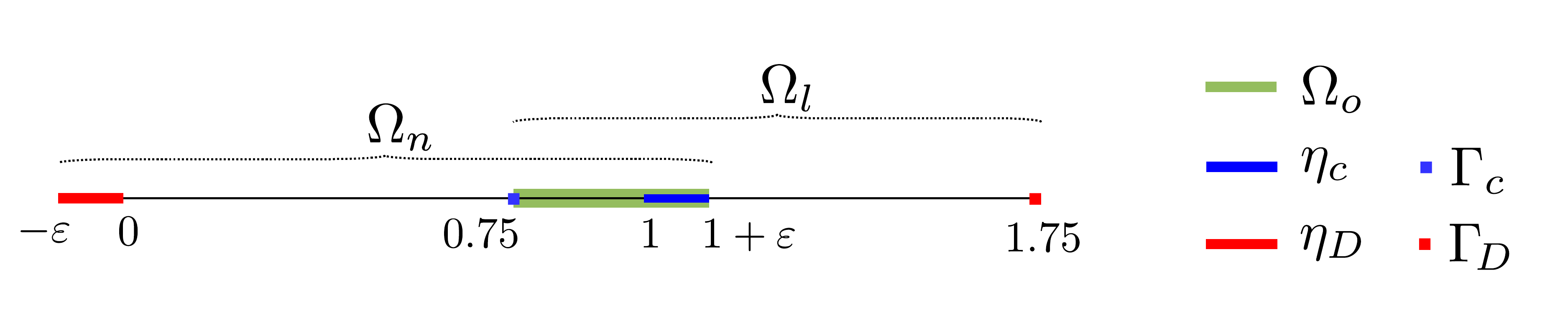}
\vspace{-2ex}
\caption{One-dimensional LtN configuration used in the numerical tests.}
\label{1D-domain}
\end{figure}

We present numerical tests with the new LtN formulation in one dimension, including a patch test, a convergence study and approximation of discontinuous solutions. Though preliminary, these results show the effectiveness of the coupling method, illustrate the theoretical results, and provide the basis for realistic simulations. 
In our examples we use an integrable kernel, $\gamma_i$, satisfying assumptions \eqref{eq:gamma-conds} and \eqref{eq:integrable-kernel} to illustrate theoretical results and a singular kernel, $\gamma_s$, often used in the literature as an approximation of a peridynamic model for nonlocal mechanics. These kernels are given by
\begin{equation}\label{test-kernel}
\gamma_i(x,y) = \frac{3}{2\varepsilon^3}\chi_{(x-\varepsilon,x+\varepsilon)}(y) 
\quad {\rm and} \quad
\gamma_s(x,y) = \frac{1}{\varepsilon^2|x-y|}\chi_{(x-\varepsilon,x+\varepsilon)}(y),
\end{equation}
respectively. Even though $\gamma_s$ does not satisfy our theoretical assumptions\footnote{The energy space associated with $\gamma_s$ is not equivalent to a Sobolev space, nevertheless it is a separable Hilbert space whose energy norm satisfies a nonlocal Poincar\'e inequality.}, these numerical results demonstrate the effectiveness of the LtN coupling for realistic, practically important, nonlocal models.
In all examples we consider the LtN problem configuration  shown in Fig.~\ref{1D-domain}, where $\omgnp=(-\varepsilon, 1+\varepsilon)$, $\omgi=(-\varepsilon,0)$, $\omgc=(1,1+\varepsilon)$, $\omgl=(0.75,1.75)$, $\gami=1.75$, $\gamc=0.75$, and $\omgb=(0.75,1+\varepsilon)$.
%
In all numerical tests $V^h_\omgi$, $V^h_\omgw$, and  $\Theta^h_n$ are discontinuous piecewise linear finite element spaces, while $H^h_\gami$ and $H^h_0$ are $C^0$ piecewise linear finite elements. We use the same grid size $h$ for the local and nonlocal finite element spaces. To solve the LtN optimization problem we apply the gradient based Quasi-Newton scheme BFGS \cite{wright:99}.

\begin{figure}
\centering
\begin{tabular}{ll}
\includegraphics[scale=.3]{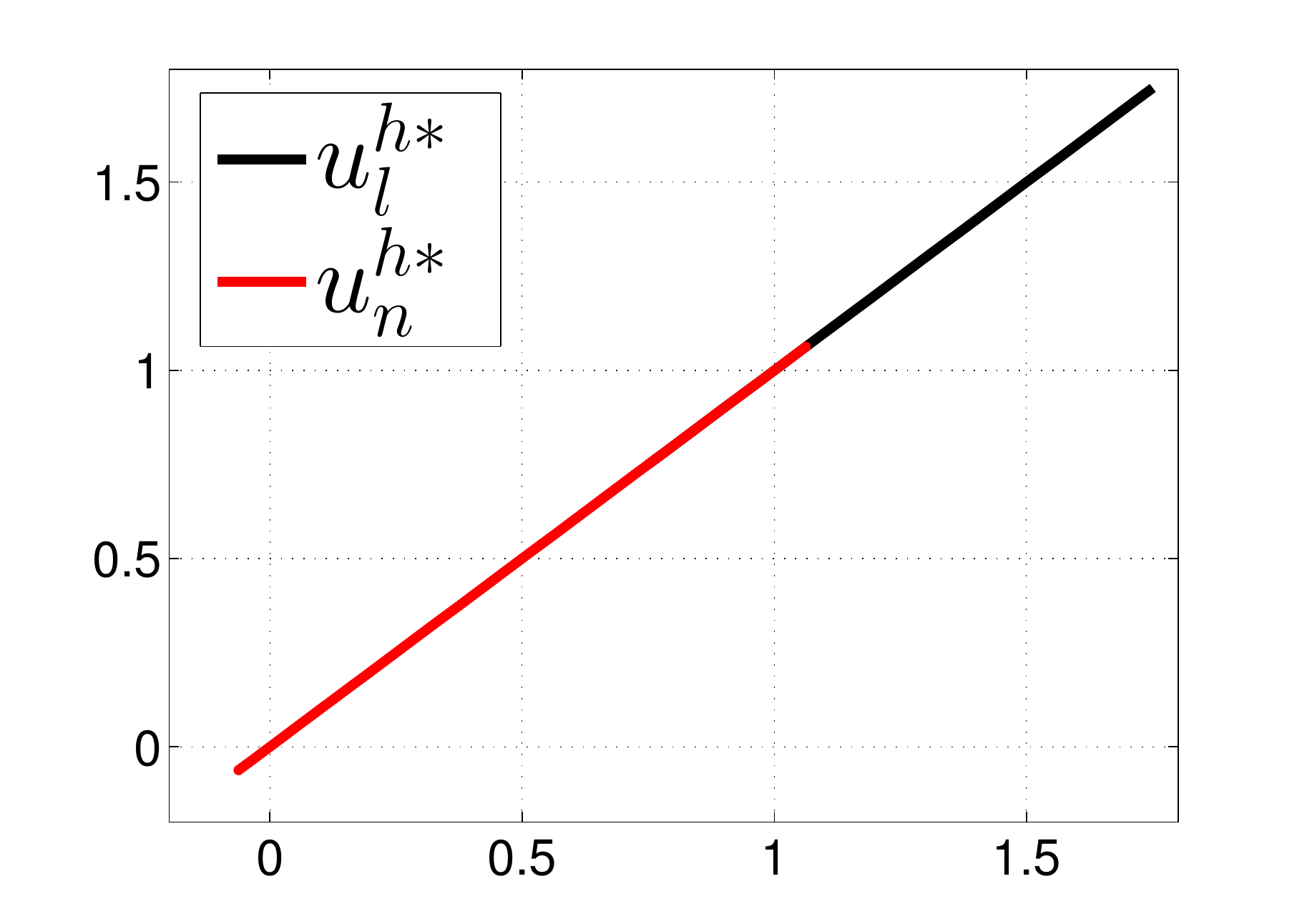} &  \hspace{-.7cm}
\includegraphics[scale=.3]{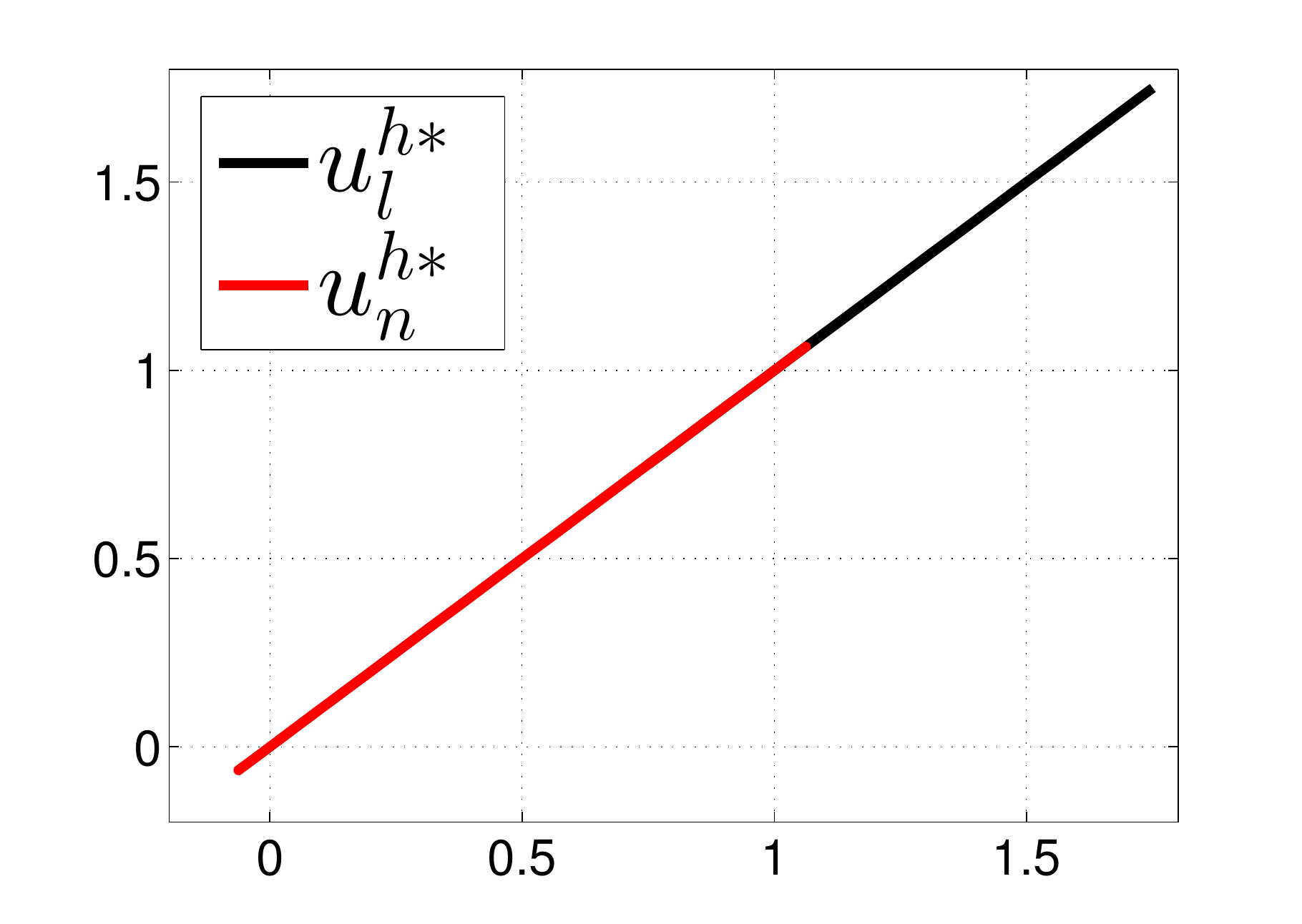}
\end{tabular}
\caption{Optimal states for the patch test with $\gamma_i$ (left) and $\gamma_s$ (right).}
\label{patch}
\end{figure}

\smallskip\paragraph{The patch test} This test uses the linear manufactured solution  $\un=\ul=x$, $\un|_\omgi=x$, $\ul(1.75)=1.75$, $f_n=f_l=0$. We expect the LtN formulation to recover this solution exactly, i.e., $u_n^{h*}\equiv\un^*$ and $u_n^{h*}\equiv\ul^*$. Figure \ref{patch} shows the optimal states $u_n^{h*}$ and $u_l^{h*}$, computed with  $\varepsilon=0.065$ and $h=2^{-7}$, for $\gamma_i$ (left) and $\gamma_s$ (right).
The LtN method recovers the exact solution to machine precision.

\begin{figure}[h!]
\centering
\begin{tabular}{ll}
\includegraphics[scale=.3]{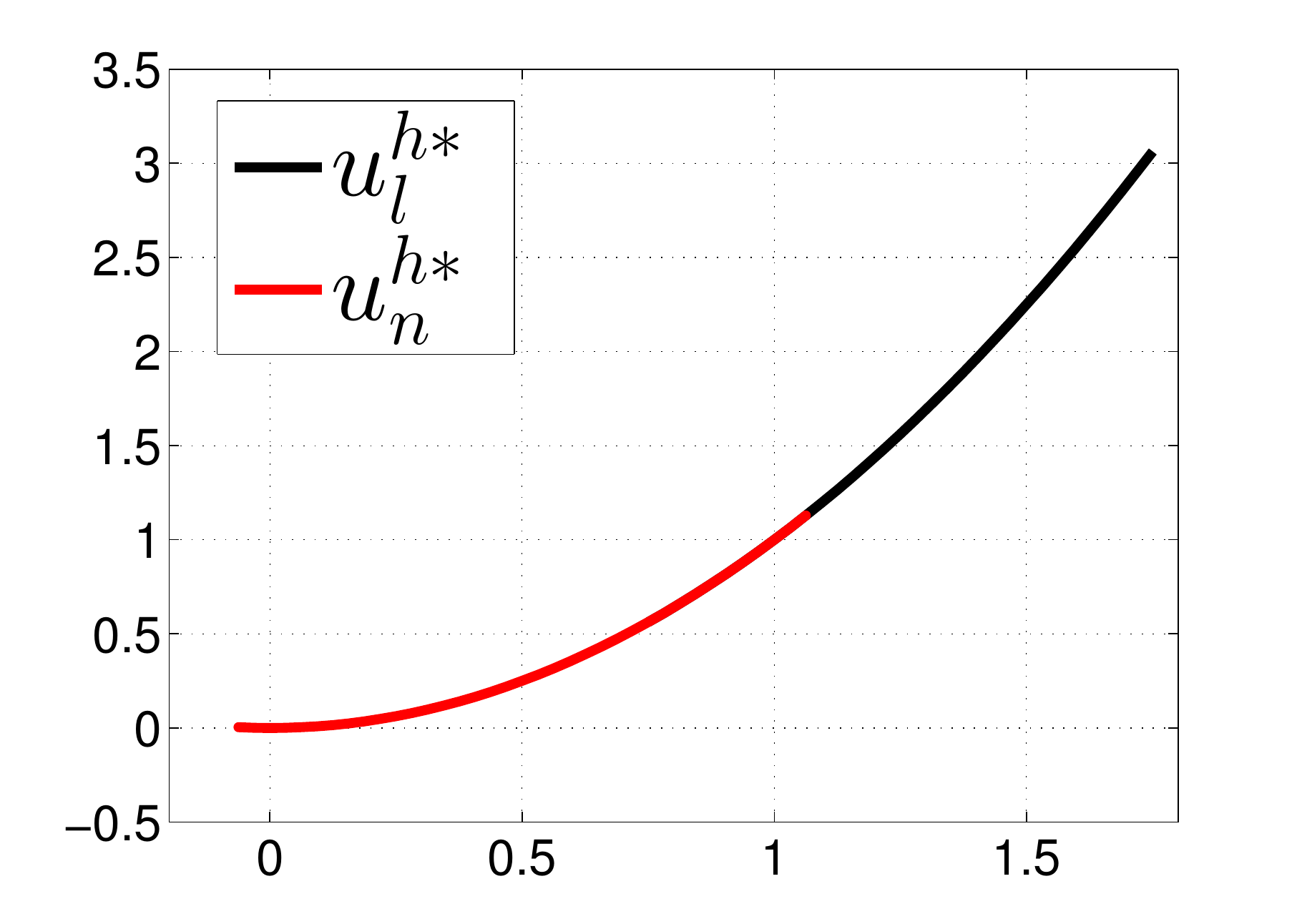} &  \hspace{-.7cm}
\includegraphics[scale=.3]{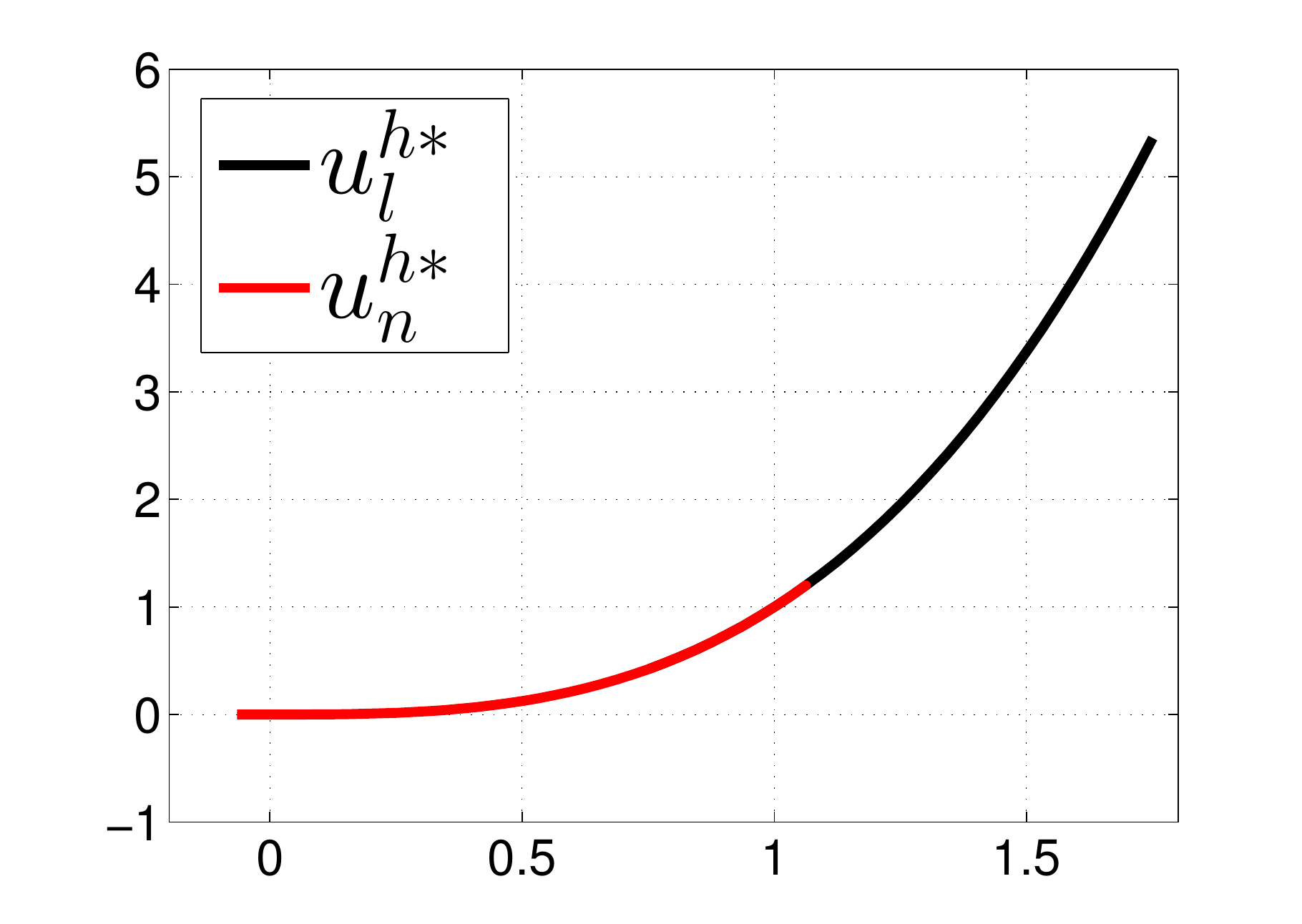}   \\\vspace{-3ex}
\includegraphics[scale=.3]{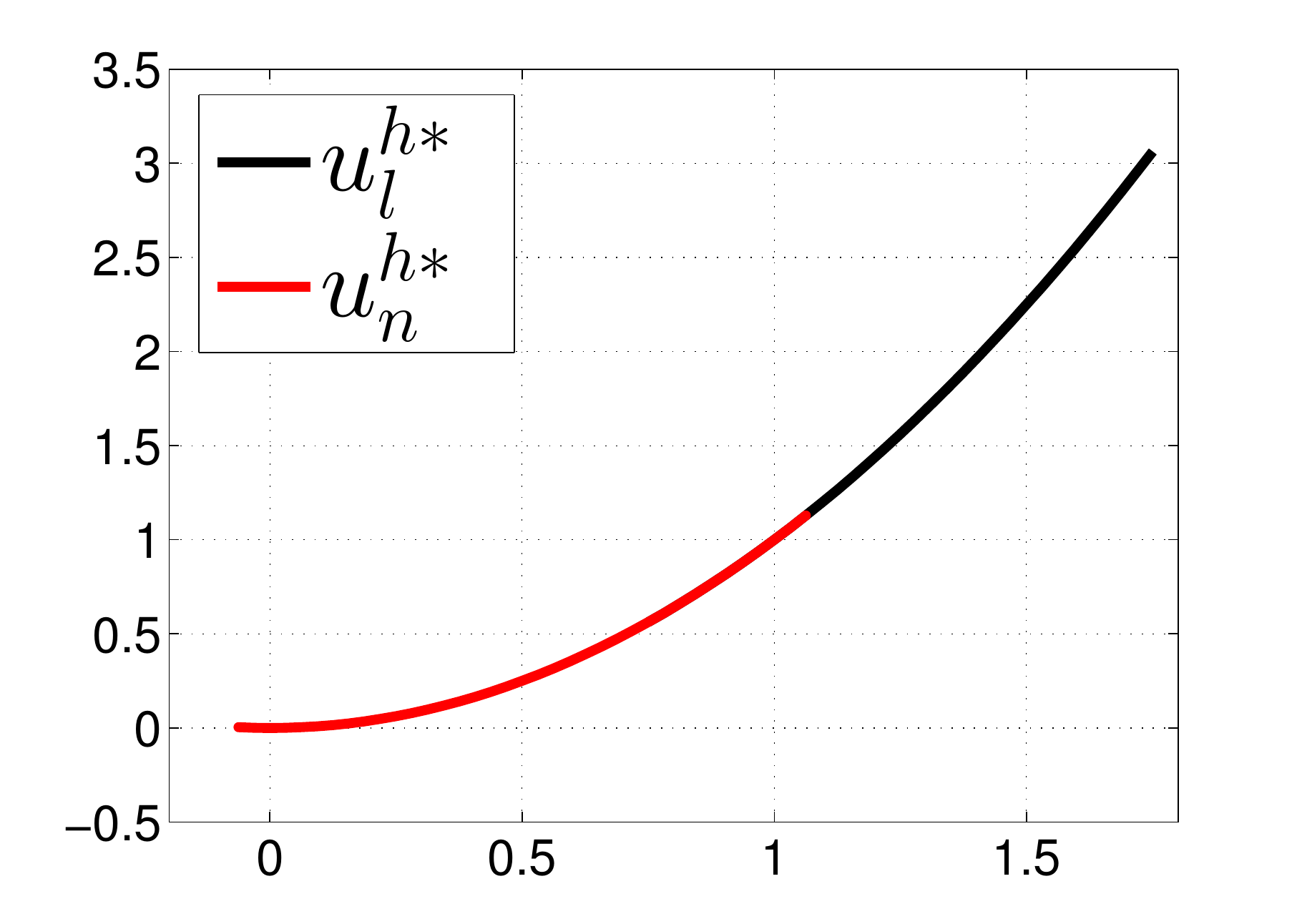} & \hspace{-.7cm}
\includegraphics[scale=.3]{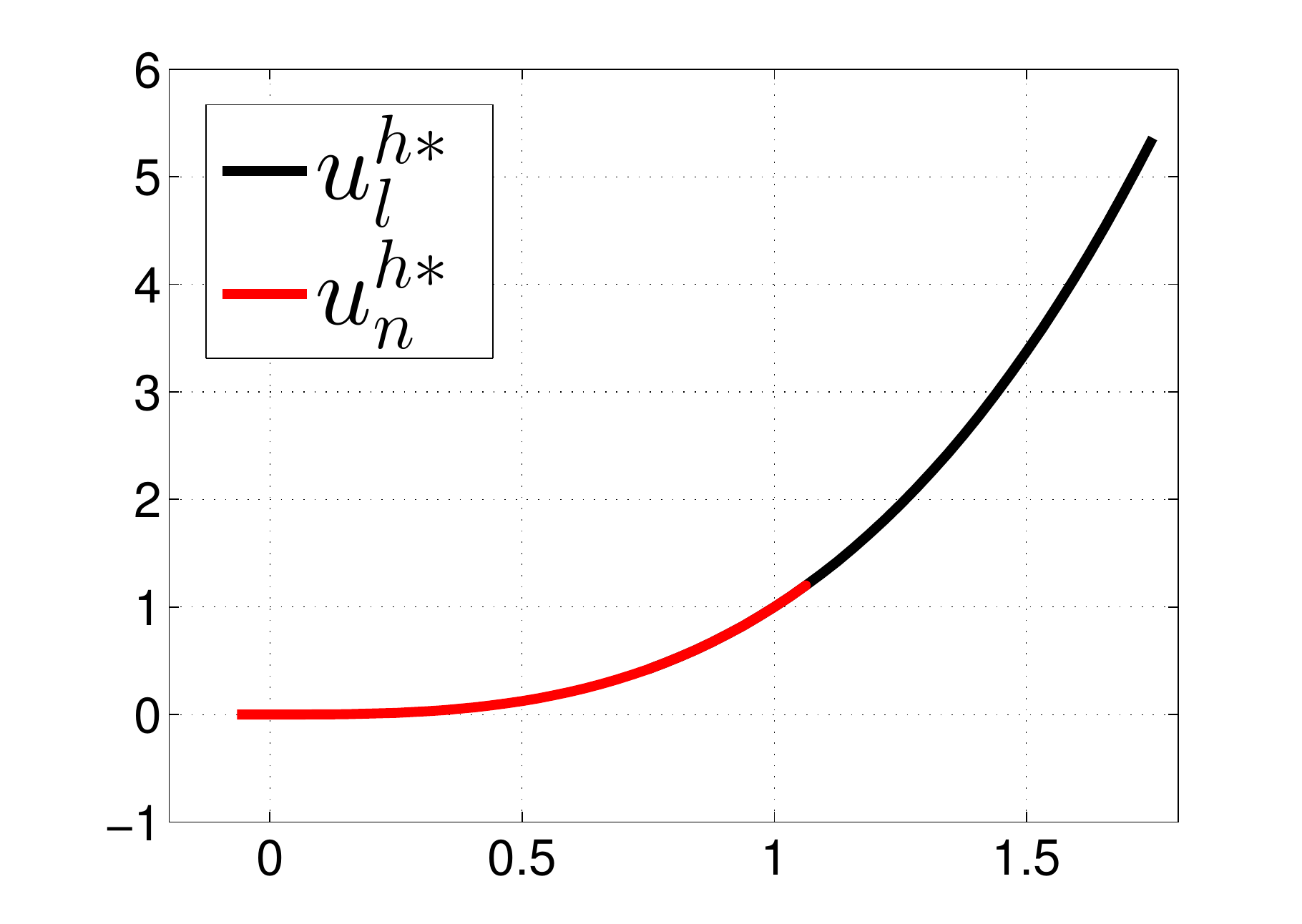}
\end{tabular}
\caption{Optimal states for \textbf{M.1} (left) and \textbf{M.2} (right) with $\gamma_i$ (top) and $\gamma_s$ (bottom).}
\label{optimal-poly}
\end{figure}

\smallskip\paragraph{Convergence tests}
 We examine the convergence of finite element approximations with respect to the grid size $h$ using the following manufactured solutions:
\begin{description}
\smallskip\item[M.1] $\un=\ul=x^2$, $\un|_\omgi=x^2$, $\ul(1.75)=1.75^2$, $f_n=f_l=-2$.
\smallskip\item[M.2] $\un=\ul=x^3$, $\un|_\omgi=x^3$, $\ul(1.75)=1.75^3$, $f_n=f_l=-6x$.
\end{description}
\smallskip Note that, for both kernels, the associated nonlocal operator is equavalent to the classical Laplacian for polynomials up to the third order.
For examples \textbf{M.1} and \textbf{M.2} we compute the convergence rates and the $L^2$ norm of the errors for the nonlocal state, $e(\un^*)$, the local state, $e(\ul^*)$, and the nonlocal control parameter, $e(\thna)$. The results are reported in Tables \ref{T:convergence2} and \ref{T:convergence3} for $\gamma_i$ and in Tables \ref{T:convergence2s} and \ref{T:convergence3s} for $\gamma_s$ in correspondence of different values of interaction radius $\varepsilon$ and grid size $h$. In Fig.~\ref{optimal-poly} we also report the optimal discrete solutions.

Results in Tables \ref{T:convergence2} and \ref{T:convergence3} show optimal convergence for state and control variables. We note that according to \cite{Du_12_SIREV} and FE convergence theory \cite{Ern_04_BOOK} this is the same rate as for the independent discretization of the nonlocal and local equations by piecewise linear elements.
\begin{remark}
The convergence analysis in Section \ref{sec:conv-anl} establishes a suboptimal convergence rate in the $L^2$ norm of the discretization error of the state variables as we lose half order of convergence. We believe that the bound in \eqref{eq:states-discr-error} is not sharp, in fact, additional numerical tests (with $h=2^{-8}, \ldots 2^{-12}$) show that there is no convergence deterioration.
\end{remark}

For the singular kernel $\gamma_s$ there are no theoretical convergence results; however, there is numerical evidence that piecewise linear approximations of \eqref{nonlocal-weak} are second-order accurate; see \cite{Chen_11_CMAME}. Our numerical experiments in Tables \ref{T:convergence2s} and \ref{T:convergence3s} show that the optimization-based LtN solution converges at the same rate.
\begin{table}[h!]
\begin{center}
\footnotesize
\begin{tabular}{| c | c | l | c | l | c | l | c |}
\hline
$\varepsilon$ &$h$ & $e(\un^*)$& rate & $e(\ul^*)$& rate & $e(\thna)$& rate \\ \hline
	& $2^{-3}$ & 2.63e-03  & -    & 2.76e-03  & -    & 5.59e-05  & -    \\  \cline{2-8}
	& $2^{-4}$ & 6.16e-04  & 2.10 & 6.74e-04  & 2.04 & 2.63e-05  & 1.09 \\  \cline{2-8}
0.010	& $2^{-5}$ & 1.40e-04  & 2.13 & 1.63e-04  & 2.05 & 1.14e-05  & 1.20 \\  \cline{2-8}  
	& $2^{-6}$ & 3.46e-05  & 2.02 & 4.04e-05  & 2.01 & 4.18e-06  & 1.47 \\  \cline{2-8}
	& $2^{-7}$ &   &  &   &  &  &  \\  \hline
	& $2^{-3}$ & 2.24e-03  & -    & 2.56e-03  & -    & 6.34e-04  & -    \\  \cline{2-8}
	& $2^{-4}$ & 7.56e-04  & 1.56 & 7.13e-04  & 1.85 & 1.78e-04  & 1.83 \\  \cline{2-8}
0.065	& $2^{-5}$ & 1.89e-04  & 2.00 & 1.78e-04  & 2.00 & 4.46e-05  & 2.00 \\  \cline{2-8}
	& $2^{-6}$ & 4.73e-05  & 2.00 & 4.46e-05  & 2.00 & 1.12e-05  & 2.00 \\  \cline{2-8}
	& $2^{-7}$ & 1.18e-05  & 2.00 & 1.11e-05  & 2.00 & 2.82e-06  & 1.99 \\  \hline
\end{tabular}
\caption{Example \textbf{M.1} with $\gamma_i$: dependence on the grid size $h$ and interaction radius $\varepsilon$ of the error.}
\label{T:convergence2}
\end{center}
\end{table}
\begin{table}[h!]
\begin{center}
\footnotesize
\begin{tabular}{ | c | c | l | c | l | c | l | c | }
\hline
$\varepsilon$ &$h$ & $e(\un^*)$& rate & $e(\ul^*)$& rate & $e(\thna)$& rate \\ \hline
	& $2^{-3}$ & 4.89e-03  & -    & 1.09e-02  & -    & 2.04e-04  & -    \\ \cline{2-8}
	& $2^{-4}$ & 1.23e-03  & 1.99 & 2.74e-03  & 2.00 & 9.63e-05  & 1.08 \\ \cline{2-8}
0.010	& $2^{-5}$ & 3.11e-04  & 1.99 & 6.86e-04  & 2.00 & 4.16e-05  & 1.21 \\ \cline{2-8}  
	& $2^{-6}$ & 7.85e-05  & 1.99 & 1.72e-04  & 2.00 & 1.45e-05  & 1.51 \\ \cline{2-8}
	& $2^{-7}$ & 1.95e-05  & 2.01 & 4.29e-05  & 2.00 & 3.16e-06  & 2.20 \\ \hline
	& $2^{-3}$ & 5.41e-03  & -    & 1.09e-02  & -    & 2.29e-03  & -    \\ \cline{2-8}
	& $2^{-4}$ & 1.34e-03  & 2.01 & 2.74e-03  & 2.00 & 5.46e-04  & 2.07 \\ \cline{2-8}
0.065	& $2^{-5}$ & 3.38e-04  & 1.99 & 6.86e-04  & 2.00 & 1.38e-04  & 1.99 \\ \cline{2-8}
	& $2^{-6}$ & 8.46e-05  & 2.00 & 1.71e-04  & 2.00 & 3.46e-05  & 1.99 \\ \cline{2-8}
	& $2^{-7}$ & 2.12e-05  & 2.00 & 4.29e-05  & 2.00 & 8.73e-06  & 1.99 \\ \hline
\end{tabular}
\caption{Example \textbf{M.2} with $\gamma_i$: dependence on the grid size $h$ and interaction radius $\varepsilon$ of the error.}
\label{T:convergence3}
\end{center}
\end{table}

\begin{table}[h!]
\begin{center}
\footnotesize
\begin{tabular}{| c | c | l | c | l | c | l | c |}
\hline
$\varepsilon$ &$h$ & $e(\un^*)$& rate & $e(\ul^*)$& rate & $e(\thna)$& rate \\ \hline
	& $2^{-3}$ & 2.67e-03  & -    & 2.78e-03  & -    & 5.79e-05  & -    \\  \cline{2-8}
	& $2^{-4}$ & 6.33e-04  & 2.08 & 6.81e-04  & 2.03 & 2.72e-05  & 1.09 \\  \cline{2-8}
0.010	& $2^{-5}$ & 1.47e-04  & 2.11 & 1.65e-04  & 2.04 & 1.19e-05  & 1.20 \\  \cline{2-8}  
	& $2^{-6}$ & 3.63e-05  & 2.01 & 4.11e-05  & 2.01 & 4.29e-06  & 1.47 \\  \cline{2-8}
	& $2^{-7}$ & 9.10e-06  & 2.00 & 1.03e-05  & 2.00 & 1.05e-06  & 2.03 \\  \hline
	& $2^{-3}$ & 2.36e-03  & -    & 2.62e-03  & -    & 6.52e-04  & -    \\  \cline{2-8}
	& $2^{-4}$ & 7.54e-04  & 1.65 & 7.12e-04  & 1.88 & 1.78e-04  & 1.87 \\  \cline{2-8}
0.065	& $2^{-5}$ & 1.88e-04  & 2.00 & 1.78e-04  & 2.00 & 4.45e-05  & 2.00 \\  \cline{2-8}
	& $2^{-6}$ & 4.67e-05  & 2.01 & 4.44e-05  & 2.00 & 1.11e-05  & 2.00 \\  \cline{2-8}
	& $2^{-7}$ & 1.14e-05  & 2.04 & 1.10e-05  & 2.01 & 2.76e-06  & 2.01 \\  \hline
\end{tabular}
\caption{Example \textbf{M.1} with $\gamma_s$: dependence on the grid size $h$ and interaction radius $\varepsilon$ of the error.}
\label{T:convergence2s}
\end{center}
\end{table}
\begin{table}[h!]
\begin{center}
\footnotesize
\begin{tabular}{ | c | c | l | c | l | c | l | c | }
\hline
$\varepsilon$ &$h$ & $e(\un^*)$& rate & $e(\ul^*)$& rate & $e(\thna)$& rate \\ \hline
	& $2^{-3}$ & 4.90e-03  & -    & 1.09e-02  & -    & 2.07e-04  & -    \\ \cline{2-8}
	& $2^{-4}$ & 1.23e-03  & 1.99 & 2.74e-03  & 2.00 & 9.68e-05  & 1.10 \\ \cline{2-8}
0.010	& $2^{-5}$ & 3.11e-04  & 1.99 & 6.86e-04  & 2.00 & 4.17e-05  & 1.21 \\ \cline{2-8}  
	& $2^{-6}$ & 7.85e-05  & 1.99 & 1.72e-04  & 2.00 & 1.46e-05  & 1.52 \\ \cline{2-8}
	& $2^{-7}$ & 1.96e-05  & 2.01 & 4.29e-05  & 2.00 & 3.17e-06  & 2.00 \\ \hline
	& $2^{-3}$ & 5.40e-03  & -    & 1.09e-02  & -    & 2.31e-03  & -    \\ \cline{2-8}
	& $2^{-4}$ & 1.34e-03  & 2.01 & 2.74e-03  & 2.00 & 5.46e-04  & 2.08 \\ \cline{2-8}
0.065	& $2^{-5}$ & 3.37e-04  & 1.99 & 6.86e-04  & 2.00 & 1.38e-04  & 2.00 \\ \cline{2-8}
	& $2^{-6}$ & 8.46e-05  & 2.00 & 1.72e-04  & 2.00 & 3.46e-05  & 1.99 \\ \cline{2-8}
	& $2^{-7}$ & 2.12e-05  & 2.00 & 4.29e-05  & 2.00 & 8.73e-06  & 1.99 \\ \hline
\end{tabular}
\caption{Example \textbf{M.2} with $\gamma_s$: dependence on the grid size $h$ and interaction radius $\varepsilon$ of the error.}
\label{T:convergence3s}
\end{center}
\end{table}

\begin{figure}
\centering
\begin{tabular}{ll}
\includegraphics[scale=.3]{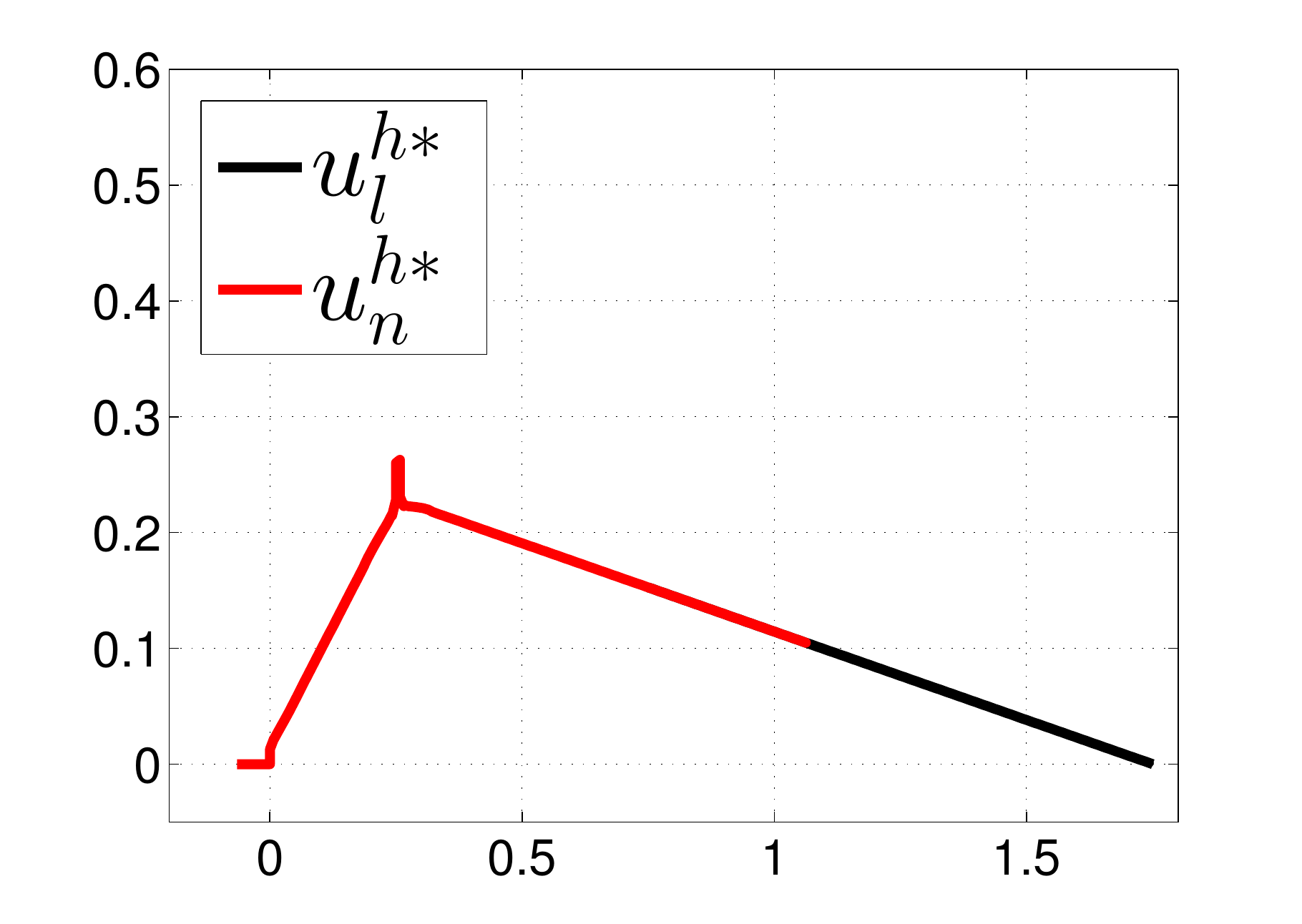} & \hspace{-.7cm}
\includegraphics[scale=.3]{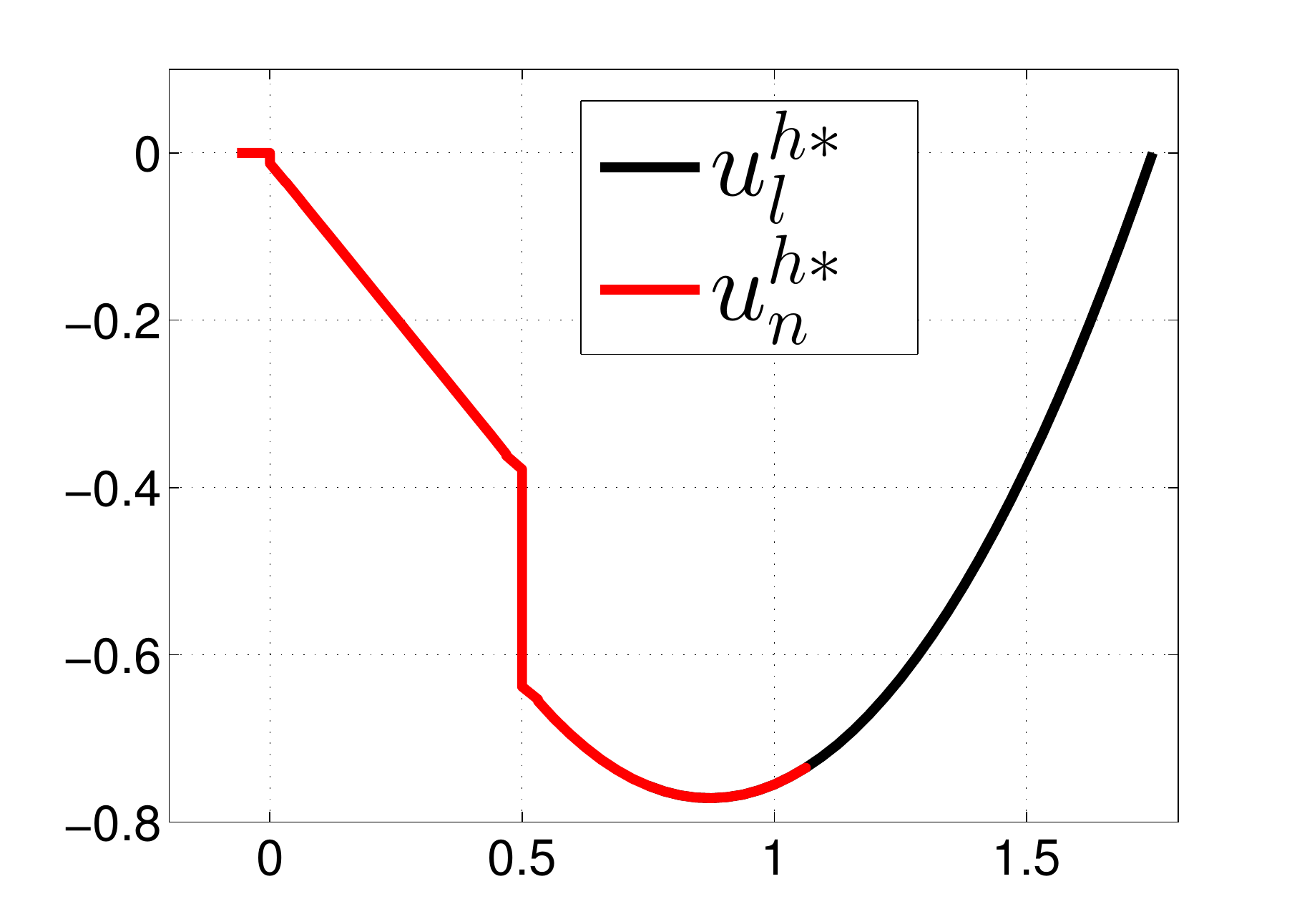} 
\end{tabular}
\vspace{-3ex}
\caption{Optimal states for examples \textbf{A.1} and \textbf{A.2}.}
\label{optimal-disc}
\end{figure}

\paragraph{Recovery of singular features} The tests in this section are motivated by nonlocal mechanics applications and demonstrate the effectiveness of the coupling method in the presence of point forces and discontinuities. We use the following two manufactured solution examples:
\begin{description}
\smallskip\item[A.1] $\un|_\omgi=0$, $\ul(1.75)=0$, $f_n=f_l=\delta(x-0.25)$, being $\delta$ the Dirac function.
\smallskip\item[A.2] $\un|_\omgi=0$, $\ul(1.75)=0$, 
\begin{displaymath}
f_n=f_l=\left\{\begin{array}{ll}
0 & x<\frac12-\varepsilon\\[4mm]
-\frac{2}{\varepsilon}\left(\frac12\varepsilon^2-\varepsilon+\frac38+\left(2\varepsilon-\frac32-\log\varepsilon\right)x\right. + & \\ 
\qquad\left.\left(\frac32+\log\varepsilon\right)x^2-\log\left(\frac12-x\right)(x^2-x)\right) & \frac12-\varepsilon\leq x<\frac12\\[4mm]
-\frac{2}{\varepsilon}\left(\frac12\varepsilon^2-\varepsilon-\frac38+\left(2\varepsilon+\frac32+\log\varepsilon\right)x\right. - & \\ 
\qquad \left.\left(\frac32+\log\varepsilon\right)x^2-\log\left(x-\frac12\right)(x^2-x)\right) & \frac12\leq x<\frac12+\varepsilon\\[4mm]
0 & x\geq \frac12 +\varepsilon.
\end{array}\right.
\end{displaymath}
\end{description}
In Fig.~\ref{optimal-disc} we report the optimal discrete solutions for $h=2^{-7}$ and $\varepsilon=0.065$. In particular, {\bf A.2} is a significant example that shows the usefulness of the coupling method in approximating the true solution with a local model where the nonlocality effects are not pronounced, i.e. the solution is smooth. 

\newpage
\section*{Acknowledgements}
This work was supported by the Sandia National Laboratories (SNL) Laboratory-directed Research and Development (LDRD) program, and the U.S. Department of Energy, Office of Science, Office of Advanced Scientific Computing Research under Award Number DE-SC-0000230927 and under the Collaboratory on Mathematics and Physics-Informed Learning Machines for Multiscale and Multiphysics Problems (PhILMs) project. Sandia National Laboratories is a multimission laboratory managed and operated by National Technology and Engineering Solutions of Sandia, LLC., a wholly owned subsidiary of Honeywell International, Inc., for the U.S. Department of Energys National Nuclear Security Administration contract number DE-NA0003525. This paper describes objective technical results and analysis. Any subjective views or opinions that might be expressed in the paper do not necessarily represent the views of the U.S. Department of Energy or the United States Government. Report Number: SAND2019-12937.

\appendix
\section{Ancillary results}\label{app:op-norm}
This appendix contains proofs of several results necessary for the estimate of $\| H \|_{**}$ and additional results necessary for the well-posedness of the discrete reduced space problem. We recall that in the following  proofs $C$ and $C_i$, $i=1,2,\ldots$, are generic positive constants.


\medskip\noindent
{\sf Proof of Lemma \ref{lemma-nonlocal-trace}.}
We need the following subspaces of $\widetilde V(\womg)$ and $V(\omgp)$
\beq\label{space-def}
\begin{array}{ll}
\widetilde V_{\womg\setminus\tau}(\womg)
&  =\{\mu\in\widetilde V(\womg),\;\mu|_{\womg\setminus\tau}=0\}\; 
    \subseteq \widetilde V(\womg), \\ 
V_{\womg\setminus\tau}(\omgp)
&=  \{w\in V({\omgp}), \; w|_{\womg\setminus\tau}=0\}\;\subseteq V(\omgp).
\end{array}
\eeq
Let $\chi_\Sigma$ be the indicator of $\Sigma$. Definition \eqref{eq:trace-norm} and the fact that $\sigma$ vanishes on $\womg\setminus\tau$ imply that
\begin{displaymath}
\|\sigma\|_{\widetilde V(\womg)}
:=\inf\limits_{v\in V({\omgp}), \,v|_\womg=\sigma} \, |||v |||_{{\omgp}}
=\inf\limits_{v\in V_{\womg\setminus\tau}(\omgp), \,v|_\womg=\sigma} \, |||v |||_{{\omgp}}\leq
\inf\limits_{v\in V_{\womg\setminus\tau}(\omgp), \,v|_\womg=\sigma} \,  |||\chi_\Sigma v |||_{{\omgp}}.
\end{displaymath}
To bound the energy norm of $\chi_\Sigma v$ note that
\begin{displaymath}
\begin{aligned}
|||\chi_\Sigma v|||^2_{\omgp} & = \int_{\omgp}\int_{\omgp}  \big(\chi_\Sigma(\xb)v(\xb)-\chi_\Sigma(\yb)v(\yb)\big)^2\gamma(\xb,\yb)\,d\yb d\xb\\
 &  = \int_\Sigma\int_\Sigma \big(v(\xb)-v(\yb)\big)^2\gamma(\xb,\yb)\,d\yb d\xb
    + \int_{\omgp\setminus\Sigma}\int_\Sigma v(\yb)^2\gamma(\xb,\yb)\,d\yb d\xb\\
 &  + \int_\Sigma\int_{\omgp\setminus\Sigma} v(\xb)^2\gamma(\xb,\yb)\,d\yb d\xb\\
 &  = |||v |||^2_\Sigma 
    + \int_\Sigma v(\xb)^2\int_{\omgp\setminus\Sigma} \gamma(\xb,\yb)\,d\yb\, d\xb
    + \int_\Sigma v(\yb)^2\int_{\omgp\setminus\Sigma} \gamma(\xb,\yb)\,d\xb\, d\yb\\[1mm]
 &  \leq |||v |||^2_\Sigma + 2\gam_1 \|v\|^2_{0,\Sigma} 
    \leq (1+ 2\gam_1 C^2_{pn}) \,|||v |||^2_\Sigma, 
\end{aligned}
\end{displaymath}
where the last two inequalities follow from \eqref{eq:gamma1} and \eqref{eq:L2equivalence} respectively. Thus,
\begin{displaymath}
\|\sigma\|_{\widetilde V(\womg)} \leq
  C\inf\limits_{v\in V_{\womg\setminus\tau}(\omgp),\,v|_\womg=\sigma}\,|||v |||_{\Sigma}
= C\inf\limits_{v\in V({\omgp}), \,v|_\womg=\sigma} \, |||v |||_{\Sigma},
\end{displaymath}
with $C^2=(1+2\gam_1 C_{pn}^2)$.
$\square$

%
\medskip\noindent
{\sf Proof of Lemma \ref{lem:closed-subsp}.}
Consider a sequence $\{\mu^k\}\subset\widetilde V_{\womg\setminus\tau}(\womg)$ such that $\mu^k\to\mu^*$ in $L^2(\womg)$. To show that $\mu^*\in \widetilde V_{\womg\setminus\tau}(\womg)$ consider the function $v^*\in L^2(\omgp)$ such that $v^* |_{\womg} = \mu^*$ and $v^* |_{\omg} = 0$.
To complete the proof it remains to show that $v^*$ has finite energy norm. Using the norm equivalence in Lemma \ref{lem:L2-equivalence}
\begin{displaymath}
|||v^* |||_\omgp \leq C^* \|v^*\|_{0,\omgp} = C^* \|\mu^*\|_{0,\womg} <\infty.
\end{displaymath}
Therefore, $v^*\in V(\omgp)$ and $v^* |_\womg = \mu^*$, hence $\mu^*\in \widetilde V_{\womg\setminus\tau}(\womg)$. Finally, let $w\in V(\omgp)$ be such that $w|_\womg=\mu$; by Lemmas \ref{lemma-nonlocal-trace} and \ref{lem:L2-equivalence}, 
$\|\mu\|_{\widetilde V(\womg)}\leq C_1 |||w|||_\womg \leq C_2\|\mu\|_{0,\womg}$.
$\square$

\medskip\noindent
{\sf Proof of Lemma \ref{lemma-strongCS}.}
We prove \eqref{c-s} by contradiction. If \eqref{c-s} does not hold then for all $0<\epsilon<1$, there exist $(\sgn^\epsilon, \sgl^\epsilon)\in \Theta_n\times\Theta_l$ such that the corresponding harmonic components $\vn(\sgn^\epsilon)=\vn^\epsilon$ and $\vl(\sgl^\epsilon)$ satisfy $\|\vn^\epsilon\|_{0,\omgb}=1$ and $\|\vl^\epsilon\|_{0,\omgb}=1$ and 
\beq\label{c-s-contradiction}
 (\vn(\sgn^\epsilon),\vl(\sgl^\epsilon))_{0,\omgb} \geq
(1-\epsilon)\|\vn(\sgn^\epsilon)\|_{0,\omgb} \, \|\vl(\sgl^\epsilon)\|_{0,\omgb}.
\eeq
Note that 
\begin{displaymath}
(1-\epsilon)\|\vn(\sgn^\epsilon)\|_{0,\omgb} \, \|\vl(\sgl^\epsilon)\|_{0,\omgb} \leq
(\vn(\sgn^\epsilon),\vl(\sgl^\epsilon))_{0,\omgb} \leq
\|\vn(\sgn^\epsilon)\|_{0,\omgb} \, \|\vl(\sgl^\epsilon)\|_{0,\omgb}.
\end{displaymath}
This implies that the sequence of inner products converges, i.e. $(\vn^\epsilon,\vl^\epsilon)_{0,\omgb}\to1$. Furthermore, since $\{\vn^\epsilon\}$ and $\{\vl^\epsilon\}$ are bounded, there exist subsequences, still denoted by $\{\vn^\epsilon\}$ and $\{\vl^\epsilon\}$, such that
\begin{equation}\label{eq:weak-epsilon-conv}
\vn^\epsilon \rightharpoonup \vn^* \quad {\rm in} \; L^2 \quad {\rm and} \quad
\vl^\epsilon \rightharpoonup \vl^* \quad {\rm in} \; H^1.
\end{equation}
Since $H^1$ is compactly embedded in $L^2$, we also have
\begin{equation}\label{eq:strong-epsilon-conv}
\vl^\epsilon \to \vl^* \quad {\rm in} \; L^2.
\end{equation}
The weak convergence of the sequences also implies that $\|\vn^*\|_{0,\omgb}\leq1$ and $\|\vl^*\|_{0,\omgb}\leq1$. Properties \eqref{eq:weak-epsilon-conv} and \eqref{eq:strong-epsilon-conv} then imply that 
\begin{displaymath}
\begin{aligned}
    \lim_{\epsilon\to1} (\vn^\epsilon,\vl^\epsilon)_{0,\omgb} 
& = \lim_{\epsilon\to1} (\vn^\epsilon-\vn^*,\vl^\epsilon)_{0,\omgb}
  + \lim_{\epsilon\to1} (\vn^*,\vl^\epsilon)_{0,\omgb} \\[2mm]
& = \lim_{\epsilon\to1} (\vn^\epsilon-\vn^*,\vl^*)_{0,\omgb}
  + \lim_{\epsilon\to1} (\vn^\epsilon-\vn^*,\vl^\epsilon-\vl^*)_{0,\omgb}
  + (\vn^*,\vl^*)_{0,\omgb} \\[2mm]
& = \lim_{\epsilon\to1} (\vn^\epsilon,\vl^\epsilon-\vl^*)_{0,\omgb}
  - \lim_{\epsilon\to1} (\vn^*,\vl^\epsilon-\vl^*)_{0,\omgb}
  + (\vn^*,\vl^*)_{0,\omgb} \\[2mm] 
& = (\vn^*,\vl^*)_{0,\omgb},
\end{aligned}
\end{displaymath}
where the last inequality follows from the strong convergence of the local component:
\begin{displaymath}
(\vn^\epsilon,\vl^\epsilon-\vl^*)_{0,\omgb} \leq
\|\vn^\epsilon\|_{0,\omgb} \|\vl^\epsilon-\vl^*\|_{0,\omgb} =
\|\vl^\epsilon-\vl^*\|_{0,\omgb} \to0 \;\; {\rm as} \;\; \epsilon\to1.
\end{displaymath}
Thus, 
\begin{displaymath}
1 =  \lim_{\epsilon\to1}(\vn^\epsilon,\vl^\epsilon)_{0,\omgb} 
  =  (\vn^*,\vl^*)_{0,\omgb}
\leq \|\vn^*\|_{0,\omgb} \|\vl^*\|_{0,\omgb} \leq 1.
\end{displaymath}
This means that $(\vn^*,\vl^*)_{0,\omgb}=\|\vn^*\|_{0,\omgb} \|\vl^*\|_{0,\omgb}$. 
This identity holds if and only if $\vl^* = \alpha \vn^*$ for some real $\alpha$. 
To complete the proof we apply the same argument as in Lemma \ref{lemma-ip}. 
On the one hand $\vn^*=0$ in $\omgi$ and so, $\vl^*=0$ in $\omgb\cap\omgi$. On the other hand, $\vl^*$ is harmonic in $\omgb$ and the identity principle implies that 
$\vl^*\equiv 0$ in $\omgb$. Thus, \eqref{c-s-contradiction} holds if and only if $\vn^*=\vl^*=0$ in $\omgb$, a contradiction.
$\square$

\medskip Next, we prove two results that are necessary for the analysis conducted in Section \ref{finite-dim}.

\smallskip
\begin{lemma}\label{lemma-discrete-SCS}
Let $\delta$ be the constant in Lemma \ref{lemma-strongCS} and let $\vnh$ and $\vlh$ be the discrete harmonic components of the discrete state variables as in \eqref{eq:discr-state-decomp}. Then, there exist $\overline h_n$ and $\overline h_l$ such that  
\begin{equation}\label{d-c-s}
|(\vnh,\vlh)_{0,\omgb}|< \delta \|\vnh\|_{0,\omgb} \, \|\vlh\|_{0,\omgb}, 
\quad \forall \;\; h_n\leq\overline h_n, \; h_l\leq\overline h_l.
\end{equation}
\end{lemma}
{\it Proof.}
In this proof we follow the same arguments of Lemma A.7 of \cite{Abdulle_15}.

Let $\{h_{nk},h_{lk}\}_{k\geq 1}$ be a sequence of mesh sizes for the nonlocal and local finite element approximations such that $\{h_{nk},h_{lk}\}\to 0$ as $k \to \infty$. Also, let $\vn^{hk}$ and $\vl^{hk}$ be the nonlocal and local discrete harmonic components corresponding to $h_{nk}$ and $h_{lk}$. It is well-known \cite{Ern_04_BOOK} that the local finite element solution converges strongly in $L^2$ to the infinite-dimensional solution $\vl$ as $k\to\infty$. Furthermore, paper \cite{Du_12_SIREV} shows that when $\gamma$ is such that \eqref{eq:integrable-kernel} holds the nonlocal finite element solution converges strongly in the energy norm to the infinite-dimensional solution $\vn$. Due to the Poincar\'e inequality such convergence implies the storng convergence in the $L^2$ norm. Thus, we have
\begin{displaymath}
\begin{aligned}
\lim_{k\to\infty}\|\vn-\vn^{hk}\|_{0,\omgb} = 0\\
\lim_{k\to\infty}\|\vl-\vl^{hk}\|_{0,\omgb} = 0.
\end{aligned}
\end{displaymath}
Since the strong convergence in $L^2$ implies the weak convergence and the convergence in norm, the following are also true
\begin{displaymath}
\begin{aligned}
\lim_{k\to\infty} (\vn^{hk},\vl^{hk})_{0,\omgb} & = & (\vn,\vl)_{0,\omgb}    \\
\lim_{k\to\infty}\|\vn^{hk}\|_{0,\omgb}      & = & \|\vn\|_{0,\omgb}  \quad\, \\
\lim_{k\to\infty}\|\vl^{hk}\|_{0,\omgb}      & = & \|\vl\|_{0,\omgb}. \quad 
\end{aligned}
\end{displaymath}
Using the strong Cauchy-Schwarz inequality \eqref{c-s}, we obtain
\begin{displaymath}
\begin{aligned}
\lim_{k\to\infty}\left|(\vn^{hk},\vl^{hk})_{0,\omgb}\right| = \left|(\vn,\vl)_{0,\omgb}\right|   
\leq \delta \, \|\vn\|_{0,\omgb} \|\vl\|_{0,\omgb} = 
     \delta \, \lim_{k\to\infty} \|\vn^{hk}\|_{0,\omgb} \|\vl^{hk}\|_{0,\omgb}. 
\end{aligned}
\end{displaymath}
Thus, $\exists \;\overline\varepsilon>0$ such that $\forall\,\varepsilon\leq\overline\varepsilon$, $\exists\; \overline h_n,\,\overline h_l$ such that $\forall\,h_n\leq \overline h_n,\;h_l\leq\overline h_l$ 
\begin{displaymath}
\left|(\vn^{hk},\vl^{hk})_{0,\omgb}\right| \leq \delta \, \|\vn^{hk}\|_{0,\omgb} \|\vl^{hk}\|_{0,\omgb}. 
\end{displaymath}
$\square$

\smallskip
\begin{lemma}\label{lemma:HilbertQh}
The space $\Theta_n^h\times\Theta_l^h$ is Hilbert with respect to the inner product 
$Q_h(\sgnh,\sglh;\mu_n^h,\mu_l^h)=(\vnh(\sgnh)-\vlh(\sglh),\vnh(\mu_n^h)-\vlh(\mu_l^h))_{0,\omgb}$.
\end{lemma}
\smallskip{\it Proof.}
By assumption $\Theta_n^h\times\Theta_l^h$ is a closed subspace $\Theta_n\times\Theta_l$. Thus, it is Hilbert with respect to $\|(\sgnh,\sglh)\|_{\Theta_n\times\Theta_l}$. Let $\|(\sgnh,\sglh)\|_{h*}=Q_h(\sgnh,\sglh;\sgnh,\sglh)$; we show that $\|(\sgnh,\sglh)\|_{\Theta_n\times\Theta_l}$ and $\|(\sgnh,\sglh)\|_{h*}$ are equivalent, i.e. there exist $K_*$ and $K^*$ such that
\begin{equation}\label{eq:energy-equiv}
K_* \|(\sgnh,\sglh)\|_{\Theta_n\times\Theta_l} \leq \|(\sgnh,\sglh)\|_{h*} \leq
K^* \|(\sgnh,\sglh)\|_{\Theta_n\times\Theta_l}.
\end{equation}
%
Using the well-posedness of the discrete problems, we have
\begin{displaymath}
\begin{aligned}
\|(\sgnh,\sglh)\|^2_{h*} & =     \|\vnh(\sgnh)-\vlh(\sglh)\|^2_{0,\omgb} 
               \leq  \|\vnh(\sgnh)\|_{0,\omgb}^2 + \|\vlh(\sglh)\|^2_{0,\omgb} \\
              & \leq C_1 \left(\|\sgnh\|^2_{0,\omgc} + \|\sglh\|^2_{\frac12,\gamc} \right).
\end{aligned}
\end{displaymath}
On the other hand, we have
\begin{displaymath}
\begin{aligned}
\|\sgnh\|^2_{0,\omgc} + \|\sglh\|^2_{\frac12,\gamc} & \leq 
\|\vnh(\sgnh)\|_{0,\omgb}^2 + C_2\|\vlh(\sglh)\|^2_{1,\omgb} & \hbox{Local trace inequality}\\
& \leq \|\vnh(\sgnh)\|_{0,\omgb}^2 + C_3\|\vlh(\sglh)\|^2_{0,\omgb} & \hbox{Local inverse inequality \cite{Ciarlet_02_BOOK}}\\
& \leq C_4\left(\|\vnh(\sgnh)\|_{0,\omgb}^2 + \|\vlh(\sglh)\|^2_{0,\omgb} \right)& \\
& \leq \frac{C_4}{1-\delta}\|\vnh(\sgnh)-\vlh(\sglh)\|^2_{0,\omgb} 
& \hbox{Strong Cauchy-Schwarz} \\
& = \frac{C_4}{1-\delta}\|(\sgnh,\sglh)\|^2_{h*}.
\end{aligned}
\end{displaymath}
Choosing $K_*=\frac{1-\delta}{C_4}$ and $K^*=C_1$, we obtain \eqref{eq:energy-equiv}. 
$\square$

\section{Notation summary}\label{app:notation}
In this appendix we report a summary of the notation we use for local and nonlocal domains. In Table \ref{T:notation} we report local entities on the left and nonlocal entities on the right (see Fig. \ref{domains} for a two-dimensional configuration). 
\begin{table}[h!]
\begin{center}
\footnotesize
\begin{tabular}{ | l | l | l | l | }
\hline
Symbol     & Definition               & Symbol     & Definition                             \\ \hline
$\Omega$   & $\omega\cup\womg$        & $\omega$   & interior of $\omgp$                    \\ 
$\Gamma$   & $\partial\omgp$          & $\eta$     & interaction domain of $\omega$         \\ \hline
$\omgl$    & local subdomain          & $\omgnp$   & nonlocal subdomain, $\omgn\cup\eta_n$  \\
           &                          & $\omgn$    & interior of $\omgnp$                   \\
$\Gamma_l$ & $\partial\omgl$          & $\eta_n$   & interaction domain of $\omgn$          \\
$\gami$    & $\Gamma\cap\Gamma_l$     & $\omgi$    & $\eta\cap\omgw$                        \\
$\gamc$    & $\Gamma_l\setminus\gami$ & $\omgc$    & $\omgw\setminus\omgi$                  \\ \hline
$\omgb$    & \multicolumn{3}{l|}{overlap domain,   $\omgnp\cap\omgl$}                       \\ \hline
\end{tabular}
\caption{Symbols used to denote local (on the left) and nonlocal (on the right) entities.}
\label{T:notation}
\end{center}
\end{table}



\end{document}